\documentclass[11pt, a4paper]{article}

\usepackage{geometry}
\usepackage{amsmath}
\numberwithin{equation}{subsection}
\usepackage{placeins}
\usepackage{titlesec}
\usepackage{multicol}
\usepackage[skip=0.1pt]{caption}
\usepackage{hyperref}

\usepackage{tikz}
\usetikzlibrary{decorations.pathreplacing,calc}

\usepackage[font=footnotesize,labelfont=bf]{caption}

\usepackage{bm}
\usepackage{natbib}
\usepackage{graphicx}
\usepackage{csquotes}
\usepackage{amsfonts}
\usepackage{xcolor}
\usepackage{tcolorbox}
\usepackage{bbm}

\usepackage{comment}
\usepackage{txfonts}

\usepackage{float}
\usepackage{algorithm}          
\usepackage{algpseudocode}      
\usepackage{lipsum} 

\titleformat*{\section}{\Large\bfseries}
\titleformat*{\subsection}{\large\bfseries}
\titleformat*{\subsubsection}{\normalsize\bfseries}

\newcommand{\beginsupplement}{%
        \setcounter{table}{0}
        \renewcommand{\thetable}{S\arabic{table}}%
        \setcounter{figure}{0}
        \renewcommand{\thefigure}{S\arabic{figure}}%
     }

\newcommand\solidrule[1][0.4cm]{\rule[0.5ex]{#1}{.4pt}}
\newcommand\dashedrule{\mbox{%
  \solidrule[1mm]\hspace{0.5mm}\solidrule[1.0mm]\hspace{0.5mm}\solidrule[1.0mm]}}

\title{Stability selection enables robust learning of partial differential equations from limited noisy data}
\author
{Suryanarayana Maddu$^{1,2,3}$, Bevan L.~Cheeseman$^{1,2,3}$, Ivo F.~Sbalzarini$^{1,2,3}$, \\Christian L. M{\"u}ller$^{4\ast}$
	\\ 
	\footnotesize{$^{1}$ \footnotesize Chair of Scientific Computing for Systems Biology, Faculty of Computer Science,}\\
	\footnotesize{ \footnotesize TU Dresden, 01069 Dresden, Germany} \\
	\footnotesize{$^{2}$ \footnotesize Center for Systems Biology Dresden, Max Planck Institute of Molecular }\\
	\footnotesize{$^{3}$\footnotesize Cluster of Excellence Physics of Life, TU Dresden, 01062 Dresden, Germany }\\
	\footnotesize{$^{4}$\footnotesize Center for Computational Mathematics, Flatiron Institute, New York, USA}
	\\
	\footnotesize{$^\ast$Correspondence to: cmueller@flatironinstitute.org}
}
\date{\vspace{-5ex}}
\begin{document}
\maketitle
\begin{abstract}
\textbf{
We present a statistical learning framework for robust identification of partial differential equations
(PDEs) from noisy spatiotemporal data. Extending previous sparse regression approaches for 
inferring PDE models from simulated data, we address key issues that have thus far limited the application 
of these methods to noisy experimental data, namely their robustness against noise and the need for manual parameter tuning. 
We address both points by proposing a stability-based model selection scheme
to determine the level of regularization required for reproducible recovery of the underlying PDE. This avoids manual parameter tuning and provides a principled way to improve the method's robustness against noise in the data. 
Our stability selection approach, termed PDE-STRIDE, can be combined with any sparsity-promoting penalized regression model and provides an interpretable criterion for model component importance. We show that in particular the combination of stability selection
with the iterative hard-thresholding algorithm from compressed sensing provides 
a fast, parameter-free, and robust computational framework for PDE inference that outperforms previous
algorithmic approaches with respect to recovery accuracy, amount of data required, and robustness to noise. 
We illustrate the performance of our approach on a wide range of noise-corrupted simulated benchmark 
problems, including 1D Burgers, 2D vorticity-transport, and 3D reaction-diffusion problems. 
We demonstrate the practical applicability of our method on real-world data by considering a purely data-driven re-evaluation of the 
advective triggering hypothesis for an embryonic polarization system in \textit{C.~elegans}. Using fluorescence microscopy images   
of \textit{C.~elegans} zygotes as input data, our framework is able to recover the PDE 
model for the regulatory reaction-diffusion-flow network of the associated proteins.
}
\end{abstract}

\section{Introduction} \label{introduction}
Predictive mathematical models, derived from first principles and symmetry 
arguments and validated in experiments, are of key importance for the scientific 
understanding of natural phenomena. While this approach has been particularly successful in 
describing spatiotemporal dynamical systems in physics and engineering, 
it has not seen the same degree of success in other scientific fields, such as neuroscience, 
biology, finance, and ecology. This is because the underlying first principles in these areas remain
largely elusive. Nevertheless, modeling in those areas has seen increasing use and relevance to help
formulate simplified mathematical equations where sufficient observational data are available for 
validation \cite{mogilner2006quantitative,sbalzarini2013modeling,tomlin2007biology,duffy2013finite, 
adomian1995solving}. In biology, modern high-throughput technologies have enabled collection of 
large-scale data sets, ranging from genomics, proteomics, and metabolomics data, 
to microscopy images and videos of cells and tissues. These data sets
are routinely used to infer parameters in hypothesized models, or to perform model selection 
among a finite number of alternative hypotheses \cite{donoho201750,Barnes2011,asmus2017p}. 
The amount and quality of biological data, as well as the performance of computing hardware and 
computational methods have now reached a level that
promises direct inference of mathematical models of biological processes from the available 
experimental data. Such data-driven approaches seem particularly valuable in cell and developmental 
biology, where first principles are hard to come by, but large-scale imaging data are available, 
along with an accepted consensus of which phenomena a model could possibly entail. 
In such scenarios, data-driven modeling approaches have the potential to uncover the unknown first principles underlying the observed biological dynamics. \\

\noindent 
Biological dynamics can be formalized at different scales, from discrete molecular processes to the 
continuum mechanics of tissues. Here, we consider the macroscopic, continuum scale where spatiotemporal
dynamics are modeled by Partial Differential Equations (PDEs) over coarse-grained state variables 
\cite{kitano2002computational,tomlin2007biology}. 
PDE models have been used to successfully address a range of biological problems 
from embryo patterning \cite{gregor2005diffusion} to modeling gene-expression 
networks \cite{chen1999modeling,ay2011mathematical} to predictive models of cell 
and tissue mechanics during growth and development \cite{prost2015active}. 
They have shown their potential to recapitulate experimental observables in cases 
where the underlying physical phenomena are known or have been postulated 
\cite{munster2019attachment}. In many biological systems, however, the governing 
PDEs are not (yet) known, which limits progress in discovering the underlying 
physical principles. Thus it is desirable to verify existing models or even discover new models by extracting governing laws directly from measured spatiotemporal data.\\

\noindent
For given observable spatiotemporal dynamics, with no governing PDE known, 
several proposals have been put forward to learn mathematically and 
physically interpretable PDE models. The earliest work in this direction 
\cite{voss1998identification} frames the problem of \enquote{PDE learning} as a
multivariate nonlinear regression problem where each component in the design matrix consists
of space and time differential operators and low-order non-linearities 
computed directly from data. Then, the alternating conditional expectation (ACE) algorithm 
\cite{Breiman1985} is used to compute both optimal element-wise non-linear transformations 
of each component and their associated coefficients. In \cite{bar1999fitting}, the problem is
formulated as a linear regression problem with a fixed pre-defined set of space and time
differential operators and polynomial transformations that are computed directly from data.
Then, backward elimination is used to identify a compact set of PDE components by minimizing
the least square error of the full model and then removing terms that reduce the fit the
least. In the statistical literature \cite{Xun2013,raissi2017machine}, 
the PDE learning problem is formulated as a Bayesian non-parametric estimation problem where 
the observed dynamics are learned via  non-parametric approximation, and a PDE representation
serves as a prior distribution to compute the posterior estimates of the PDE coefficients. 
Recent influential work revived the idea of jointly learning the structure {\em 
and} the coefficients of PDE models from data in discrete space and 
time using sparse regression ~\cite{brunton2016discovering,rudy2017data, schaeffer2017learning}. Approaches such as SINDy \cite{brunton2016discovering} and PDE-FIND 
\cite{rudy2017data} compute a large pre-assembled dictionary of possible PDE terms from data 
and identify the most promising components via penalized linear regression formulations. 
For instance, PDE-FIND is able to learn different types of PDEs from simulated spatiotemporal
data, including Burgers, Kuramato-Sivishinksy, reaction-diffusion, and Navier-Stokes 
equations. PDE-FIND's performance was evaluated on noise-free simulated data as well as data 
with up to 1\% additive noise and showed a critical dependency on the proper  
setting of the tuning parameters which are typically unknown in practice. 
Recent approaches attempt to alleviate this dependence by using Bayesian sparse 
regression schemes for model uncertainty quantification \cite{Zhang2018e} or 
information criteria for tuning parameter selection \cite{Mangan2019}.
There is also a growing body of work that considers deep neural networks for PDE learning 
\cite{long2017pde,raissi2018hidden,Raissi2019}. For instance, a deep feed forward network 
formulation \cite{long2017pde}, PDE-NET, directly learns a computable discretized form of the
underlying governing PDEs for forecasting ~\cite{long2017pde, long2018pde}. 
PDE-NET exploits the connection between differential operators and 
orders-of-sum rules of convolution filters \cite{dong2017image} in order to 
constrain network layers to learn valid discretized differential operators. The forecasting 
capability of this approach was numerically demonstrated for predefined linear differential 
operator templates. A compact and interpretable symbolic identification of the PDE 
structure is, however, not available with this approach.
\\

\noindent 
Here, we ask the question whether and how it is possible to extend the class of sparse regression inference methods to work on real, limited, and noisy experimental data.
As a first step, we present a statistical learning framework, PDE-STRIDE (STability-based 
Robust IDEntification of PDEs), to robustly infer PDE models from noisy spatiotemporal data
without requiring manual tuning of learning parameters, such as regularization constants. 
PDE-STRIDE is based on the statistical principle of stability selection 
\cite{meinshausen2010stability,shah2013variable}, which provides an interpretable criterion
for any component's inclusion in the learned PDE in a data-driven manner. Stability 
selection can be combined with any sparsity-promoting regression method, including 
LASSO~\cite{tibshirani1996regression, meinshausen2010stability}, iterative hard 
thresholding (IHT)~\cite{blumensath2008iterative}, Hard Thresholding Pursuit 
(HTP)~\cite{foucart2011hard}, or Sequential Thresholding Ridge Regression 
(STRidge)~\cite{rudy2017data}. PDE-STRIDE therefore provides a drop-in solution to render 
existing inference tools more robust, while reducing the need for parameter tuning. In the benchmarks presented 
herein, the combination of stability selection with de-biased iterative hard thresholding (IHT-d)
empirically showed the best performance and highest consistency w.r.t.~perturbations 
of the dictionary matrix and sampling of the data.\\

\noindent This paper is organized as follows: Section~\ref{methods} provides the 
mathematical formulation of the sparse regression problem and discusses how the design 
matrix is assembled. We also review the concepts of regularization paths and stability 
selection and discuss how they are combined in the proposed method. The numerical results 
in Section~\ref{results} highlight the performance and robustness of
the PDE-STRIDE for recovering different PDEs from noise-corrupted 
simulated data. We also perform an achievability analysis of PDE-STRIDE+IHT-d inference scheme for consistency and convergence of recovery probability with increasing
sample size. Section~\ref{result_5} demonstrates that the robustness of the 
proposed method is sufficient for real-world applications. We consider learning a PDE 
model from noisy biological microscopy images of membrane protein dynamics in a 
\textit{C.~elegans} zygote. Section~\ref{conclusion} provides a summary of our 
results and highlights future challenges for data-driven PDE learning.

\begin{figure}[h]
\begin{minipage}[b]{1.0\textwidth}
\centering
\includegraphics[width=0.75\linewidth]{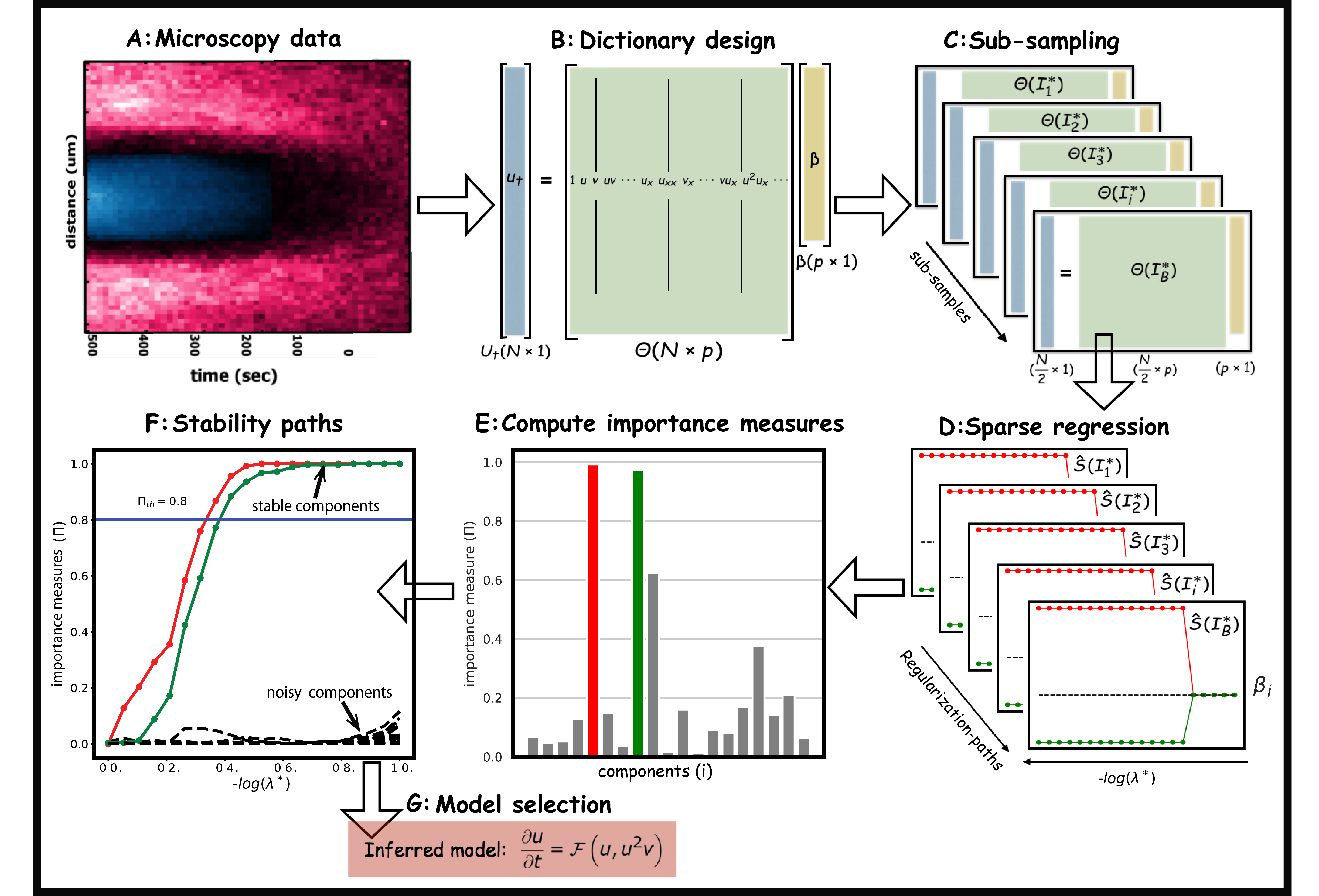}
\end{minipage}
\vspace{0.0em}
\captionsetup{justification=raggedright,margin=0.5cm}
\caption{\textbf{Enabling data-driven mathematical model discovery through stability selection:}
We outline the necessary steps in our method for learning PDE models from spatiotemporal data. 
\textbf{A}: Extract spatiotemporal profiles from microscopy videos of the chosen state
variables. Data courtesy of Grill Lab, MPI-CBG/TU Dresden \cite{gross2019guiding}. \textbf{B}:
Compile the design matrix $\Theta$ and the measurement vector $U_t$ from the data. \textbf{C}:
Construct multiple linear systems of reduced size through random sub-sampling of the rows of 
the design matrix $\Theta$ and $U_t$. \textbf{D}: Solve and record the sparse/penalized regression 
solutions independently for each sub-sample along the $\lambda$ paths. \textbf{E}: Compute the 
importance measure $\Pi$ for each component. The histogram shows the importance measure $\Pi$ for all components at a particular value of $\lambda$. \textbf{F}: Construct the stability plot by 
aggregating the importance measures along the $\lambda$ path, leading to separation of the 
noise variables (dashed black) from the stable components (colored). Identify the most stable components by 
thresholding $\Pi > 0.8$.
\textbf{G}: Build the PDE model from the identified components.}
\label{schematic}
\end{figure}
\section{Problem Formulation and Optimization} \label{methods}
We outline the problem formulation underlying the data-driven PDE inference considered 
here. We review important sparse regression techniques and introduce the concept of 
stability selection used in PDE-STRIDE. 

\subsection{Problem formulation for PDE learning} \label{formulation}
We propose a framework for stable estimation of the structure and parameters of the 
governing equations of continuous dynamical systems from discrete spatiotemporal 
measurements or observations. Specifically, we consider PDEs for the multidimensional state
variable $u \in \mathbb{R}^d$ of the form shown in Eq. (\ref{general_parametric_form}),
composed of polynomial non-linearities (e.g., $u^2, u^3$), spatial derivatives (e.g.,
$u_{x},u_{xx}$), and the parametric dependence modeled through $\Xi \in \mathbb{R}^{p}$. \begin{equation}\label{general_parametric_form}
    \frac{\partial u }{\partial t} = \mathcal{F}\left( [u, u^2,u^{3}, u_{xx}, uu_{x},.....], x,t,\Xi \right).
\end{equation}
Here, $\mathcal{F}(\cdot)$ is the function map that models the spatiotemporal non-linear 
dynamics of the system. For simplicity, we limit ourselves to forms of the function map 
$\mathcal{F}(\cdot)$ that can be written as the linear combinations of polynomial 
non-linearities, spatial derivatives, and combinations of both. For instance, for 
a one-dimensional ($d=1$) state variable $u$, the function map can take the form:
\begin{equation}\label{explicit_govern}
\frac{\partial u}{\partial t} = \underbrace{\xi_0 + \xi_1  u + \xi_{2}\frac{\partial u}{\partial x} + \xi_3 u \frac{\partial u}{\partial x} + \xi_4 u^2 + ..... + \xi_{k}  u^{3} \frac{\partial^2 u}{\partial x^2 } + ....,}_{\mathcal{F}(\cdot)}
\end{equation}
\noindent where $\xi_k$ are the coefficients of the components of the PDE for $k\ge 0$. 
The continuous PDE of the form described in Eq.~(\ref{explicit_govern}), with appropriate
coefficients $\xi_k$, holds true for all continuous space and time points $(x,t)$ in the
domain of the model. Numerical solutions of the PDE try to satisfy the equality relation 
in Eq. (\ref{explicit_govern}) for reconstituting the non-linear dynamics of a dynamical
system at discrete space and time points $(x_m, t_n)$.
We assume that we have access to $N$ noisy observational data $\tilde u_{m,n}$ of the 
state variable $u$ over the discrete space and time points. The measurement errors are 
independent and identically distributed and are assumed to follow a normal distribution 
with mean zero and variance $\sigma^2$. \\

\noindent We follow the formulation put forward in \cite{bar1999fitting,rudy2017data, 
schaeffer2017learning} and construct a large dictionary of PDE components using  
discrete approximations of the dynamics from data. 
For instance, for the one-dimensional example in Eq. (\ref{explicit_govern}),
the discrete approximation with $p$ PDE terms can be written in vectorized form as a linear regression problem: 
\begin{equation}
\label{eq:linReg}
\underbrace{\begin{bmatrix}
           \vert \\
           u_t \\
           \vert
\end{bmatrix}}_\text{$U_t$} 
=
\underbrace{\begin{bmatrix}
          \:\: \vert \qquad \vert \qquad \vert \qquad \vert \qquad \vert \qquad \vert \quad\\
          1 \qquad \!\! u \quad uu_{x}  \quad \:\:...\:\: \:  u^3 u_{xx} \quad ..\\
          \:\: \vert \qquad \vert \qquad \vert \qquad \vert \qquad \vert \qquad \vert \quad
\end{bmatrix}}_\text{$\Theta$}
\times \,  \xi.
\end{equation}

\noindent Here, the left-hand side $U_t \in \mathbb{R}^{N}$ is a discrete 
approximation to the time derivatives of $u$ at the discretization points and 
represents the response or outcome vector in the linear regression framework. 
Each column of the dictionary or design matrix $\Theta \in \mathbb{R}^{N\times 
p}$ represents the discrete approximation of one PDE component, i.e., one of the
terms in Eq. (\ref{explicit_govern}), at $N$ discrete points in space and time
$(x, t)_n ,\, n=1,\ldots , N$. Each column is interpreted as a potential
predictor of the response vector $U_t$. The vector ${\xi} = \left[ \xi_0,\,
\xi_1,\, \ldots \, \xi_{p-1} \right]^{\top} \in \mathbb{R}^{p}$ is the vector of
unknown PDE coefficients, i.e., the pre-factors of the terms in Eq.
(\ref{explicit_govern}).\\

\noindent Both $U_t$ and $\Theta$ need to be constructed from numerical 
approximations of the temporal and spatial derivatives of the observed state 
variables. There is a rich literature in numerical analysis on this topic (see, 
e.g.,\cite{chartrand2011numerical, stickel2010data}). 
Here, we approximate the time derivatives by 
first-order forward finite differences from $\tilde u$ after initial denoising of
the data. Similarly, the spatial derivatives are computed by applying the
second-order central finite differences. Details about the denoising are given in Section \ref{subsub:compData} and the Supplemental Material.\\ 

\noindent Given the general linear regression \emph{ansatz} in 
Eq.~\ref{eq:linReg} we formulate the data-driven PDE inference problem as a 
regularized optimization problem of the form: 
\begin{equation}\label{general_opt}
\hat{\xi}^{\lambda} = \arg \min_{\xi} \Big( h(\xi) + \lambda g(\xi) \Big) \,,
\end{equation}
where $\hat{\xi}^{\lambda}$ is the minimizer of the objective function, $h(\cdot)$ a 
smooth convex data fitting function, $g(\cdot)$ a regularization or penalty 
function, and $\lambda \ge 0$ is a scalar tuning parameter that balances data 
fitting and regularization. The function $g(\cdot)$ is not necessarily convex or
differentiable. We follow previous work \cite{bar1999fitting,rudy2017data, 
schaeffer2017learning} and consider the standard least-squares objective 
\begin{equation}\label{square_loss}
 h(\xi) = \frac{1}{2}\Vert U_t - \Theta \xi \Vert_{2}^{2} \,.
\end{equation}
The choice of the penalty function $g(\cdot)$ influences the properties of the 
coefficient estimates $\hat{\xi}^{\lambda}$. Following 
\cite{brunton2016discovering,rudy2017data,schaeffer2017learning}, 
we seek to identify a small set of PDE components among the $p$ possible 
components that accurately predict the time evolution of the state variables. 
This implies that we want to identify a sparse coefficient vector $\hat 
\xi^{\lambda}$, thus resulting in an \textit{interpretable} PDE model. This can 
be achieved through sparsity-promoting penalty functions $g(\cdot)$. We next 
consider different choices for $g(\cdot)$ that enforce sparsity in the 
coefficient vector and review algorithms that solve the associated optimization 
problems.

 

\subsection{Sparse optimization for PDE learning}
\label{sparse_promoters}


\noindent 
The least-squares loss in Eq.(\ref{square_loss}) can be combined with different sparsity-inducing penalty functions $g(\cdot)$. The prototypical example is the $\ell_1$-norm
$g(\cdot) = \Vert \cdot \Vert_1$ leading to the LASSO formulation of sparse linear regression \cite{tibshirani1996regression}:
\begin{equation}\label{lasso}
\hat{\xi}^{\lambda} = \arg \min_{\xi} \Bigg( \underbrace{ \frac{1}{2} \Vert U_t -
\Theta \xi \Vert_{2}^{2}}_{h(\cdot)} +  \underbrace{\lambda \Vert \xi 
\Vert_1}_{g(\cdot)} \Bigg) \, ,
\end{equation}
The LASSO objective comprises a convex smooth loss and a convex non-smooth 
regularizer. For this class of problems, efficient optimization algorithms exist
that can exploit the properties of the functions and come with convergence 
guarantees. 
Important examples include coordinate-descent algorithms 
\cite{Wu2008b,friedman2010regularization} and proximal algorithms, including the 
Douglas-Rachford algorithm \cite{eckstein1992douglas} and the projected (or proximal) 
gradient method, also known as the Forward-Backward algorithm 
\cite{combettes2011proximal}. In signal processing, the latter schemes are termed 
iterative shrinkage-thresholding (ISTA) algorithms (see \cite{Beck2009} and references 
therein) which can be extended to non-convex penalties. 
Although the LASSO has been previously used for PDE learning in 
\cite{schaeffer2017learning}, the statistical performance of the LASSO estimates are 
known to deteriorate if certain conditions on the design matrix are not met. For 
example, the studies in \cite{meinshausen2006high, zhao2006model} provide sufficient 
and necessary conditions, called the \textit{irrepresentable conditions}, for 
consistent variable selection using LASSO, essentially excluding too strong 
correlations of the predictors in the design matrix. The conditions are, however, 
difficult to check in practice, as they require knowledge of the true components of the
model. One way to relax these conditions is via randomization. The Randomized LASSO \cite{meinshausen2010stability} considers the following objective:
\begin{equation}\label{random_lasso}
\hat{\xi}^{\lambda} = \arg \min_{\xi} \Bigg( \underbrace{ \frac{1}{2} \Vert U_t - \Theta \xi \Vert_{2}^{2}}_{h(\cdot)} +  \underbrace{\lambda \sum_{k=1}^{p}   \frac{\vert \xi_k \vert}{W_k}}_{g(\cdot)} \Bigg) \, ,
\end{equation}
where each $W_k$ is an {\it i.i.d.}~random variable uniformly distributed over 
$[\alpha,1]$ with $\alpha \in (0,1]$. For $\alpha = 1$, the Randomized LASSO reduced to
the standard LASSO. The Randomized LASSO has been shown to successfully overcome
limitations of the LASSO to handle correlated components in the dictionary 
\cite{meinshausen2010stability} while simultaneously preserving the overall 
convexity of objective function. As part of our PDE-STRIDE framework, we will evaluate
the performance of the Randomized Lasso in the context of PDE learning using cyclical
coordinate descent \cite{friedman2010regularization}.\\

\noindent
The sparsity-promoting property of the (weighted) $\ell_1$-norm comes at the expense 
of considerable bias in the estimation of the non-zero coefficients \cite{zhao2006model},
thus leading to reduced variable selection performance in practice. 
This drawback can be alleviated by using non-convex penalty functions
\cite{Fan2001,Zhang2010c}, allowing near-optimal variable selection performance at the
cost of needing to solve a non-convex optimization problem. For instance, using 
the $\ell_0$-``norm" (which counts the number of non-zero elements of a vector) as
regularizers $g(\cdot) = \Vert \cdot \Vert_0$ leads to the NP-hard problem:
\begin{equation}\label{hard_problem}
    \hat{\xi}^{\lambda} = \arg \min_{\xi} \bigg( \underbrace{ \frac{1}{2} \Vert U_t - \Theta \xi \Vert_{2}^{2}}_{h(\cdot)} +  \underbrace{ \lambda \Vert \xi \Vert_0 }_{g(\cdot)} \bigg) \, .
\end{equation}
This formulation has found widespread applications in compressive sensing 
and signal processing. Algorithms that deliver approximate solutions to the objective 
function in Eq.~\ref{hard_problem} include greedy optimization strategies 
\cite{tropp2004greed}, Compressed Sampling Matching Pursuit (CoSaMP)
\cite{needell2009cosamp}, and subspace pursuit \cite{dai2009subspace}.
We here consider the Iterative Hard Thresholding (IHT) algorithm 
\cite{blumensath2008iterative, blumensath2009iterative}, which belongs to the 
class of ISTA algorithms. Given the design matrix $\Theta$ and the measurement vector 
$U_t$, IHT computes sparse solutions $\hat \xi$ by applying a non-linear shrinkage (thresholding) operator to gradient descent steps in an iterative manner.
One step in the iterative scheme reads:
\begin{equation}\label{IHT}
\xi^{n+1} = H_{\lambda} \left ( \xi^{n} + \Theta^{T}(U_t - \Theta \xi^{n}) \right) = T \xi^{n}, \quad \quad H_{\lambda} (x)= 
\begin{cases}
    0,& x \leq \sqrt{\lambda}\\
    x,              & \text{otherwise}
\end{cases}.
\end{equation}
The operator $H_{\lambda}(\xi)$ is the non-linear hard-thresholding operator. Convergence of the 
above IHT algorithm is guaranteed iff $ \Vert U_t - \Theta \xi^{n+1} \Vert_{2}^{2} < (1-c)\Vert 
U_t - \Theta \xi^{n} \Vert_{2}^{2} $ is true in each iteration for some constant $0 < c < 1$.  
Specifically, under the condition
that $\Vert \Theta \Vert_2 < 1$, the IHT algorithm is guaranteed to not increase the cost function
in Eq.~(\ref{hard_problem}) (\textit{Lemma 1} in  \cite{blumensath2008iterative}). The above IHT
algorithm can also be viewed as a thresholded version of the classical Landweber iteration
\cite{herrity2006sparse}. The fixed points $\xi^{*}$ of $\xi^{*} = T \xi^{*}$ for the non-linear
operator $T$ in Eq.~(\ref{IHT}) are local minima of the cost function in Eq.~(\ref{hard_problem})
(\textit{Lemma 3} in  \cite{blumensath2008iterative}). Under the same condition on the design
matrix, i.e.~$ \Vert \Theta \Vert_2 < 1 $, the optimal solution of the cost function in
Eq.~(\ref{hard_problem}) thus belongs to the set of fixed points of the IHT algorithm
(\textit{Theorem 2} in \cite{blumensath2008iterative} and \textit{Theorem 12} in
\cite{tropp2006just}). Although the IHT algorithm comes with theoretical convergence 
guarantees, the resulting fixed points are not necessarily sparse 
\cite{blumensath2008iterative}. \\


\noindent Here, we propose modification of the IHT algorithm that will prove to be 
particularly suited for solving PDE learning problems with PDE-STRIDE. Following a 
proposal in \cite{foucart2011hard} for the Hard Thresholding Pursuit (HTP) algorithm, 
we equip the standard IHT algorithm with an additional debiasing step. At each 
iteration, we solve a least-squares problem restricted to the support $S^{n+1} = \{ 
k:\xi^{n+1} \neq 0 \}$ obtained from the $n^{th}$ IHT iteration. We refer to this form 
of IHT as Iterative Hard Thresholding with debiasing (IHT-d). In this two-stage 
algorithm, the standard IHT step serves to extract the explanatory variables, while the 
debiasing step approximately debiases (or refits) the coefficients restricted to the 
currently active support \cite{figueiredo2007gradient}. Rather than solving the 
least-squares problem to optimality, we use gradient descent steps until a loose upper 
bound on the least-squares refit is satisfied. This prevents over-fitting by attributing low 
confidence to large supports and reduces computational overhead. The complete IHT-d procedure
is detailed in Algorithm 1 in the Supplementary material. In PDE-STRIDE, we will compare
IHT-d with a heuristic iterative algorithm, Sequential Thresholding of Ridge regression 
(STRidge), that also uses $\ell_{0}$ penalization and is available in PDE-FIND 
\cite{rudy2017data}.\\

\subsection{Stability selection}\label{stability_selection}
The practical performance of the sparse optimization techniques for PDE learning critically 
hinges on the proper selection of the regularization parameter $\lambda$ that balances model
fit and complexity of the coefficient vector. In model discovery tasks on real experimental 
data, a wrong choice of the regularization parameter could result in incorrect PDE model 
selection even if true model discovery would have been, in principle, achievable. In 
statistics, a large number of tuning parameter selection criteria are available, ranging 
from cross-validation approaches \cite{Kohavi1995} to information criteria 
\cite{Schwarz1978}, or formulations that allow joint learning of model coefficients and 
tuning parameters \cite{Lederer2015,Bien2016o}. Here, we advocate stability-based model 
selection \cite{meinshausen2010stability} for robust PDE learning. 
The statistical principle of stability \cite{Yu2013b} has been put forward 
as one of the pillars of modern data science and statistics and provides an intuitive
approach to model selection \cite{meinshausen2010stability,shah2013variable,Liu2010}.\\ 

\noindent In the context of sparse regression, stability selection \cite{meinshausen2010stability} 
proceeds as follows (see also Figure \ref{schematic} for an illustration). 
Given a design matrix $\Theta$ and the measurement vector $U_t$, we generate random subsample 
indices $I_i^* \subset \{1,\ldots,\}, \, i=1,\ldots ,B$ of equal size $\vert I_i^* \vert = N/2$ 
and produce reduced sub-designs $\Theta [I_{i}^{*}] \in \mathbb{R}^{\frac{N}{2} \times p}$ and 
$U_t [I_{i}^{*}] \in \mathbb{R}^{\frac{N}{2}}$ by choosing rows according to the index set $I_i^*$. 
For each of the resulting $B$ subproblems, we apply a sparse regression technique
and systematically record the recovered support $\hat{S}^{\lambda}[I_i^*], \,  i=1,\ldots ,B$ 
as a function of $\lambda$ over an regularization path $\Lambda = [ \lambda_{\textrm{max}}\lambda_{\textrm{min}} ]$. 
The values of $\lambda_{\textrm{max}}$ are data dependent and are 
easily computable for generalized linear models with convex penalties 
\cite{friedman2010regularization}. Similarly, the critical parameter $\lambda_{\textrm{max}}$ for
the non-convex problem in Eq.(\ref{hard_problem}) can be evaluated from optimality conditions 
(\textit{Theorem 12} in \cite{tropp2006just} and \textit{Theorem 1} in 
\cite{blumensath2008iterative}). The lower bound $\lambda_{\textrm{min}}$ of the regularization 
path is set to $\lambda_{\textrm{min}} = \epsilon \lambda_{\textrm{max}}$ with default value 
$\epsilon = 0.1$.  The $\lambda$-dependent stability (or importance) measure for each
coefficient $\xi_k$ is then computed as:
\begin{equation}\label{importance}
    \hat{\Pi}_{k}^{\lambda} = \mathbb{P}\left(k \in \hat{S}^{\lambda}\right) \approx \frac{1}{B} \sum_{i=1}^{B} \mathbbm{1}(k \in \hat{S}^{\lambda}[I_{i}^{*}] ),
\end{equation}
\noindent where $I_{1}^{*},\ldots ,I_{B}^{*}$ indicates the independent random sub-samples. 
The stability $\hat{\Pi}_{k}^{\lambda}$ of each coefficient can be plotted across the 
$\lambda$-path, resulting in a component stability profile (see Figure \ref{schematic}F for an 
illustration). This visualization provides an intuitive overview of the importance of the different 
model components. Different from the original stability selection proposal 
\cite{meinshausen2010stability}, we define the stable components of the model as follows: 
\begin{equation}\label{stable_components}
    \hat{S}_\text{stable} = \{ k: \hat{\Pi}_{k}^{\lambda_{\textrm{min}}}  \geq \pi_{th}\} 
\end{equation}
\noindent Here, $\pi_{th}$ denotes the critical stability threshold parameter and can be set to 
$\pi_{th} \in [0.7, 0.9]$ \cite{meinshausen2010stability}. The default setting is $\pi_{th}=0.8$.
In an exploratory data analysis mode, the threshold $\pi_{th}$ can also be set through visual 
inspection of the stability plots, thereby allowing the principled exploration of several 
alternative PDE models. The stability measures $\hat{\Pi}^\lambda_{k}$ also provide an interpretable criterion for a model component's  importance, thereby guiding the user to build 
the right model with high probability. As we will show in the numerical experiments, 
stability selection ensures robustness against varying dictionary size, different types of data 
sampling, noise in the data, and variability of the sub-optimal solutions when non-convex 
penalties are used. All of these properties are critical for consistent and reproducible model 
learning in real-world applications. Under certain conditions, stability selection can also 
provides an upper bound on the expected number of false positives. Such guarantees are not 
generally assured by any sparsity-promoting regression method in 
isolation~\cite{shah2013variable}. For instance, stability selection combined with randomized 
LASSO (Eq.~(\ref{random_lasso}) with $\alpha < 0.5$) is consistent for variable selection even 
when the irrepresentable condition is violated~\cite{meinshausen2010stability}.

\section{Numerical experiments on simulation data} \label{results}
We present numerical experiments in order to benchmark the performance and robustness of PDE-STRIDE combined with different $\ell_{0}/\ell_{1}$ sparsity-promoting regression methods to infer PDEs from spatiotemporal data. In order to provide comparisons and benchmarks, we first use simulated data obtained by numerically solving known ground-truth PDEs, before applying our method to a real-world data set from biology. The benchmark experiments on simulation data are presented in four sub-sections that demonstrate different aspects of the inference framework: Sub-section~\ref{result_1} demonstrates the use of  different sparsity-promoting regression methods in our framework in a simple 1D Burgers problem. Sub-section~\ref{result_2} then compares their performance in order to choose the best regression method, IHT-d. In sub-section \ref{result_3}, stability selection is combined  with IHT-d to recover 2D vorticity-transport and 3D reaction-diffusion PDEs from limited, noisy  simulation data. Sub-section~\ref{result_4} reports  achievability results to quantify the robustness of stability selection to perturbations in dictionary size, sample size, and noise levels. \\

\begin{tcolorbox}
  \footnotesize{
\textbf{STability-based 
Robust IDEntification of PDEs (PDE-STRIDE)}\\
Given the noise-corrupted data $\tilde u$ and a choice of regression method, e.g., (randomized) LASSO, IHT, HTP, IHT-d, STRidge.
\begin{enumerate}
    \item Apply any required de-noising method on the noisy data and compute the spatial derivatives and non-linearities to construct the design matrix $\Theta \in \mathbb{R}^{N \times p}$ and the time-derivatives vector $U_t \in \mathbb{R}^{N \times 1}$ for suitable sample size and dictionary size, $N$ and $p$, respectively.
    \item Build the sub-samples $\Theta[I_{i}^{*}] \in \mathbb{R}^{\frac{N}{2} \times p}$ and $U_t[I_{i}^{*}]$, for $i = {1,2,...,B}$, by uniformly randomly sub-sampling of rows from the design matrix $\Theta$ and the corresponding rows from $U_t$. For every sub-sample $I_i^{*}$, standardize the sub-design matrix $\Theta[I_i^{*}]$ such that  $\sum_{j=1}^{\frac{N}{2}} \theta_{jk} = 0$ and $\frac{1}{N} \sum_{j=1}^{\frac{N}{2}} \theta_{jk}^{2} = 1$, for $k = 1,2,...,p$. Here, $\theta_{jk}$ is the element in row $j$ and column $k$ of the matrix $\Theta[I_i^{*}]$. The corresponding measurement vector $U_t[I_i^{*}]$ is centered to zero mean.
    \item Apply the sparsity-promoting regression method independently to each sub-sample $\Theta[I_i^{*}], U_t[I_i^{*}]$ to construct the $\lambda$ paths for $M$ values of $\lambda$ as discussed in section \ref{stability_selection}.
    \item Compute the importance measures $\hat{\Pi}_{k}^{\lambda}$ of all dictionary components $\xi _k$ along the discretized $\lambda$ path, as discussed in section \ref{stability_selection}. Select the stable support set $\hat{S}_\text{stable}$ by applying the threshold $\pi_{th}$ to all $\hat{\Pi}_k$. Solve a linear least-squares problem restricted to the support $\hat{S}_\text{stable}$ to identify the coefficients of the learned model.
\end{enumerate}
}
\end{tcolorbox}

\subsubsection*{Adding noise to the simulation data}
Let $u \in \mathbb{R}^{N}$ be the vector of clean simulation data sampled in both space and time. This vector is corrupted with additive Gaussian noise to 
\begin{equation*}
    \tilde u = u + \varepsilon,
\end{equation*}
such that $\varepsilon = \sigma \cdot \mathcal{N}\left(0,  \textrm{std}(u)\right) $ is the additive Gaussian noise with an empirical standard deviation of the entries in the vector $u$,  and $\sigma$ is the level of Gaussian noise added. 

\subsubsection*{Computing the data vector}
\label{subsub:compData}
The data vector $U_t\in \mathbb{R}^{N}$ contains numerical approximations to the time derivatives of the state variable $u$ at different points in space and time. We compute these by first-order forward finite differences (i.e., the explicit Euler scheme)
from $\tilde u$ after initial de-noising of the data. Similarly, the spatial derivatives are computed by applying the second-order central finite differences directly on the de-noised data. For de-noising we use truncated single value decomposition (SVD) with a cut-off at the elbow of the singular values curve, as shown in Supplementary Figures~\ref{Burgers_SVD} and  \ref{GSRD_SVD}.

\subsubsection*{Fixing the parameters for stability selection}
We propose that PDE-STRIDE combined with IHT-d provides a parameter-free PDE learning method. Therefore, all stability selection parameters described in Section \ref{stability_selection} are fixed throughout our numerical experiments. The choice of these statistical parameters is well-discussed in the 
literature \cite{meinshausen2010stability, buhlmann2014high, friedman2010regularization}. We thus fix: the repetition number $B = 250$, regularization path parameter $\epsilon=0.1$, $\lambda$-path size $M = 20$, and the importance threshold $\pi^{th} = 0.8$. Using these standard choices, the methods works robustly across all tests presented hereafter, and is parameter-free in that sense. In both stability and regularization plots, we show the normalized value of regularization constant $\lambda^{*} = \lambda/\lambda_{\textrm{max}}$. Although, the stable component set $\hat{S}$ is evaluated at $\lambda_{\textrm{min}}$ as in Eq.(\ref{stable_components}), the entire stability profile of each component from $\lambda_{\textrm{max}}$ to $\lambda_{\textrm{min}}$ is shown in all our stability plots. This way, we get additional graphical insight into how each component evolves along the $\lambda-$path.

\subsection{Case study with 1D Burgers equation and different sparsity-promoters} \label{result_1}

\noindent We show that stability selection can be combined with any sparsity-promoting penalized regression to learn PDE components from noisy and limited spatiotemporal data. We use simulated data of the 1D Burgers equation 
\vspace{-0.5em}
\begin{equation}\label{Burgers}
    \frac{\partial u}{\partial t} + u \frac{\partial u}{\partial x} = \frac{\partial^2 u }{\partial x^2}
\end{equation}
with identical boundary and initial conditions as used in \cite{rudy2017data} 
to provide fair comparison between methods: periodic boundaries in space and the following Gaussian initial condition:
\begin{equation*}
    u(x,0) = e^{ \left( -(x+2)^2 \right) }, \quad x \in [-8,8]
\end{equation*}
The simulation domain $ [ -8, 8] $ is divided uniformly into 256 Cartesian grid points in space and 1000 time points. 
The numerical solution is visualized in space-time in Figure \ref{Burgers_results_p19}. The numerical solution was obtained using parabolic method based on finite differences and time-stepping using explicit Euler method with step size $dt = 0.01$.\\

\begin{figure}[h]
\begin{minipage}[b]{1.0\textwidth}
\centering
\includegraphics[width=0.75\linewidth]{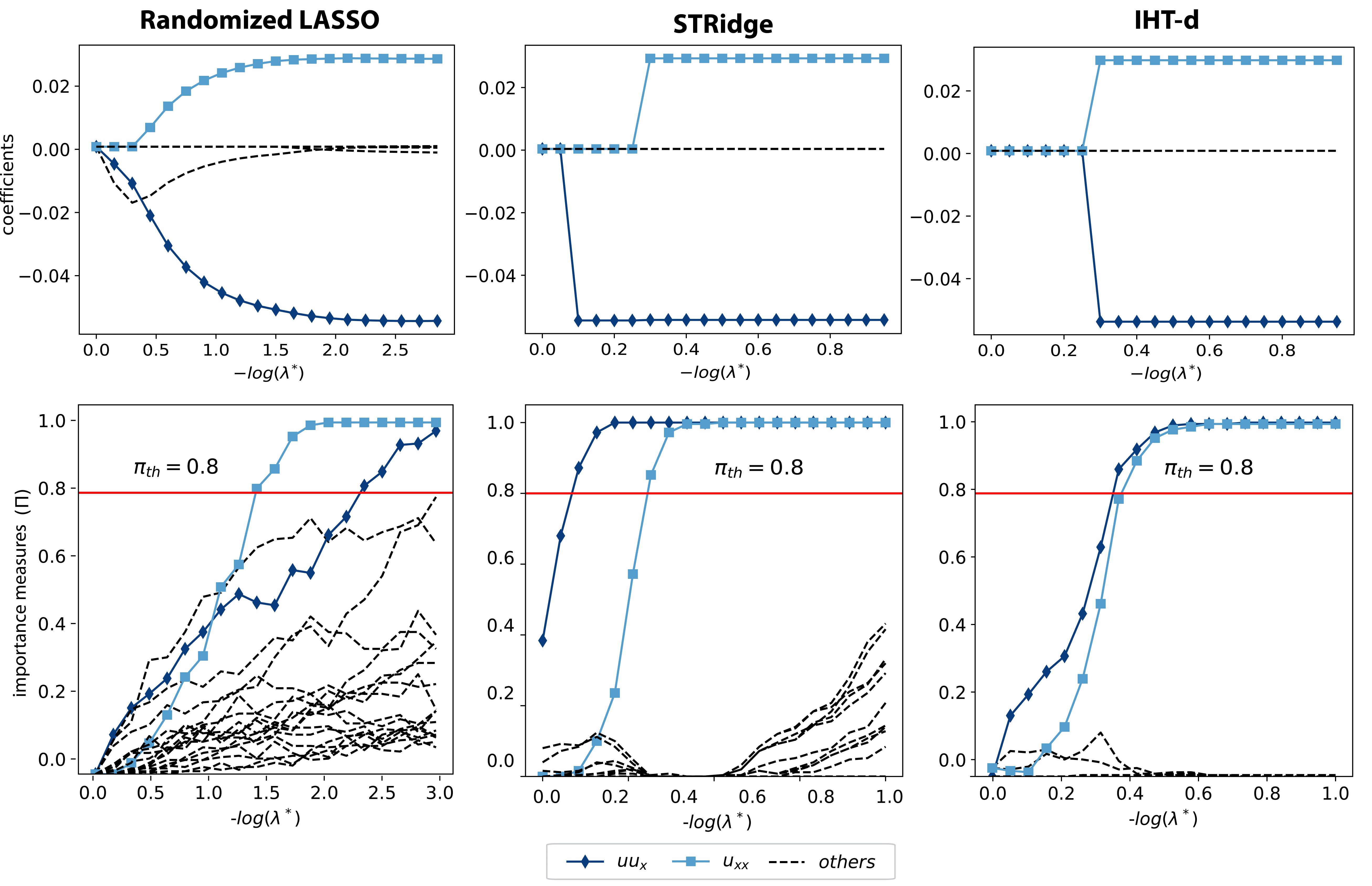}
\end{minipage}
\captionsetup{justification=raggedright,margin=0.5cm}
\caption{ \textbf{Model selection with PDE-STRIDE for the 1D Burgers equation:} The top row shows regularization paths (see Section \ref{stability_selection}) for three sparsity-promoting regression techniques: randomized LASSO, STRidge, and IHT-d all for the  same design ($N = 200$, $p =19$). 
 The inset at the bottom shows the colors correspondence with the dictionary. The ridge parameter $\lambda_{R}$ for STRidge is fixed to $\lambda_R = 10^{-5}$ \cite{rudy2017data}. The value of $\alpha$ for the randomized LASSO is set to $0.2$.
In all three cases, the standard threshold $\pi_{th}=0.8$ (red solid line, (\textcolor{red}{\solidrule})) correctly identifies the true components. NOTE: The $\epsilon$ is set to $0.001$ for the randomized LASSO case for demonstrating stability selection.}
\label{Result_1}
\end{figure}

\noindent We test the combinations of stability selection with the three sparsity-promoting regression techniques described
in Section \ref{sparse_promoters}: randomized LASSO, STRidge, and IHT-d. The top row of Figure \ref{Result_1} shows the regularization paths
for randomized LASSO, STRidge, and IHT-d. The bottom row of Figure 2 shows the corresponding stability profiles for each component in the dictionary. The colored solid lines correspond to
the advection and diffusion terms of Eq.\ref{Burgers}.
\noindent
Thresholding the importance measure at $\Pi > \pi_{th}$ = 0:8, sparsity-promoting regression methods are able to identify the correct components in the stability plots (solid colored lines) from the noise variables (dashed black lines). 

\subsection{Comparison between sparsity-promoting techniques} \label{result_2}
\begin{figure}[h]
\begin{minipage}[b]{1.0\textwidth}
\centering
\includegraphics[width=0.7\linewidth]{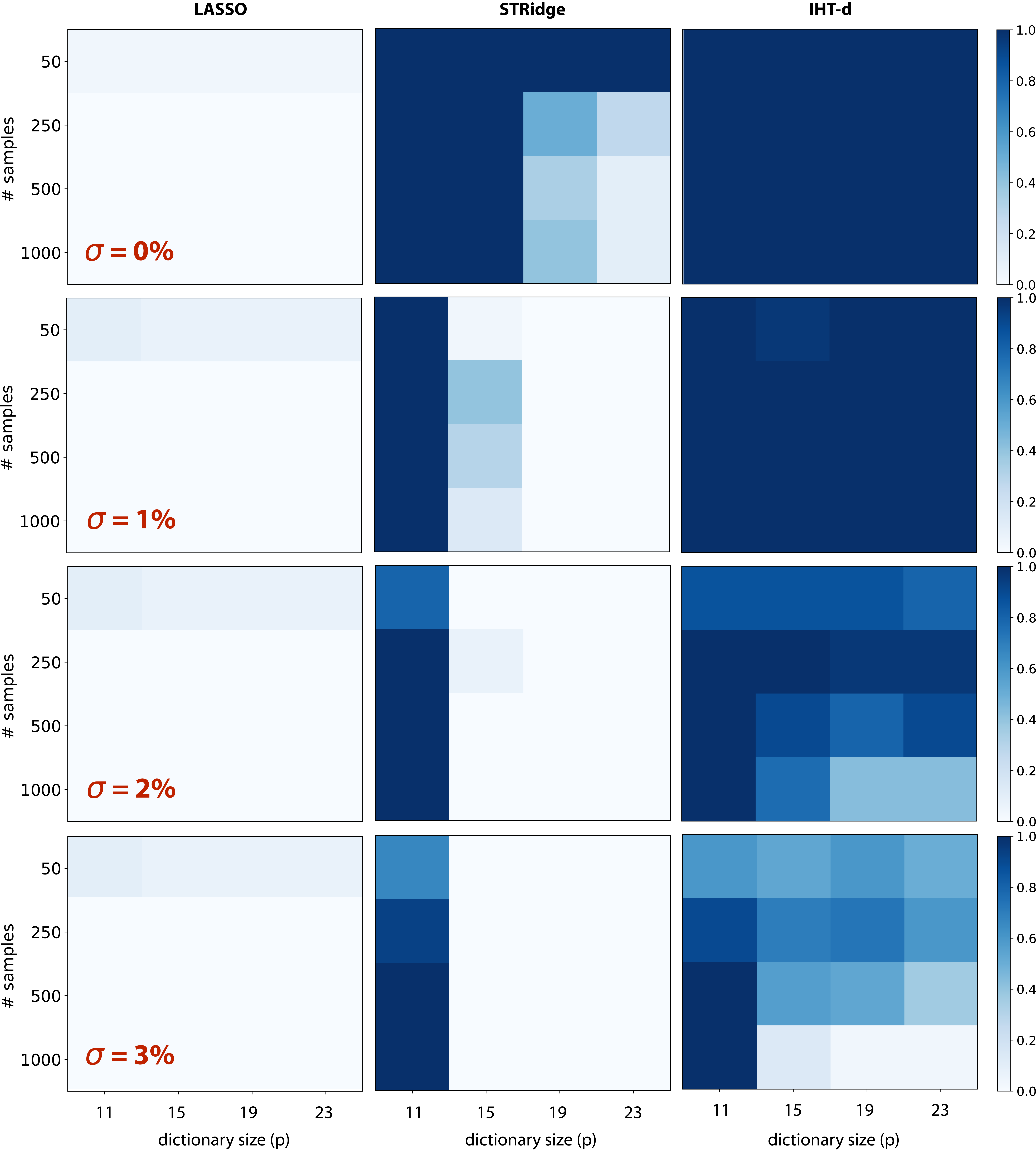}
\end{minipage}
\\
\captionsetup{justification=raggedright,margin=0.5cm}
\caption{ \textbf{Comparison between different sparsity promoters for the 1D Burgers equation:}  Each colored square corresponds to a design $(N, p, \sigma)$ with certain sample size $N$, disctionary size $p$, and nose level $\sigma$. Color indicates the success frequency over 30 independent repetitions with uniformly random data sub-samples. ``Success'' is defined as the existence of {\em a} $\lambda$ for which the correct PDE is recovered from the data. The columns compare the three sparsity-promoting techniques: LASSO, STRidge, and IHT-d (left to right), as labelled at the top. }
\label{Result_bestsubset}
\end{figure}

Although stability selection can be used in conjunction with any $\ell_1$ or $\ell_0$ sparsity-promoting regression method, the question arises whether a particular choice of regression algorithm is particularly well suited in the frame of PDE learning. We therefore perform a systematic comparison between LASSO,
STRidge, and IHT-d for recovering the 1D Burgers equation under perturbations to the sample size $N$, the dictionary size $p$, and the noise level $\sigma$. An experiment for a particular triple $(N, p, \sigma)$ is considered as success if there exists a $\lambda \in \Lambda$ (see section \ref{stability_selection}) for which the true PDE components are recovered. In Figure \ref{Result_bestsubset}, the success frequencies over 30 independent repetitions with uniformly random data sub-samples are shown.
\\

\noindent A first observation from Figure \ref{Result_bestsubset} is that $\ell_{0}$  solutions (here with STRidge and IHT-d) are better than the relaxed $\ell_{1}$ solutions (here with LASSO). We also observe that IHT-d performs better than LASSO and STRidge for  large dictionary sizes $p$, high noise, and small sample sizes.
Large dictionaries with higher-order derivatives computed from discrete data cause grouping (correlations between variables), for which LASSO tends to select one variable from each group, ignoring the others~\cite{wang2011random}. Thus, LASSO fails to identify the true support consistently. STRidge shows good recovery for large dictionary sizes $p$ with clean data, but it breaks down in the presence of noise on the data. Of all three methods, IHT-d shows the best robustness to both noise and changes in the design. We note a decrease in inference power with increasing sample size $N$, especially for large $p$ and high noise levels. This again can be attributed to correlations and groupings in the dictionary, which become more prominent with increasing sample size $N$. \\

\noindent Based on these results, we use IHT-d as the sparsity-promoting regression method in conjunction with PDE-STRIDE for model selection in the remaining sections. 

\subsection{Stability-based model inference } \label{result_3}
We present benchmark results of PDE-STRIDE for PDE recovery with IHT-d as the sparse regression method. This combination of methods is used to recover PDEs from limited noisy data obtained by numerical solution of the 1D Burgers, 2D vorticity-transport, and 3D Gray-Scott equations. Once the support $\hat{S}$ of the PDE model is learned by PDE-STRIDE with IHT-d, the actual coefficient values of the non-zero components are determined by solving the linear least-squares problem restricted to the recovered support $\hat{S}$. However, more sophisticated methods could be used for parameter estimation for a known structure of the PDE like in \cite{raissi2017machine,Xun2013} from limited noisy data. But, this is beyond the scope of this paper given that in all cases considered the sample size $N$ significantly exceeds the cardinality of the recovered support ($N \gg \vert \hat{S} \vert$) for which LLS fit provide good estimates of the PDE coefficients.

\subsubsection{1D Burgers equation}
\begin{figure}[h]
\begin{minipage}[b]{1.0\textwidth}
\centering
\includegraphics[width=0.75\linewidth]{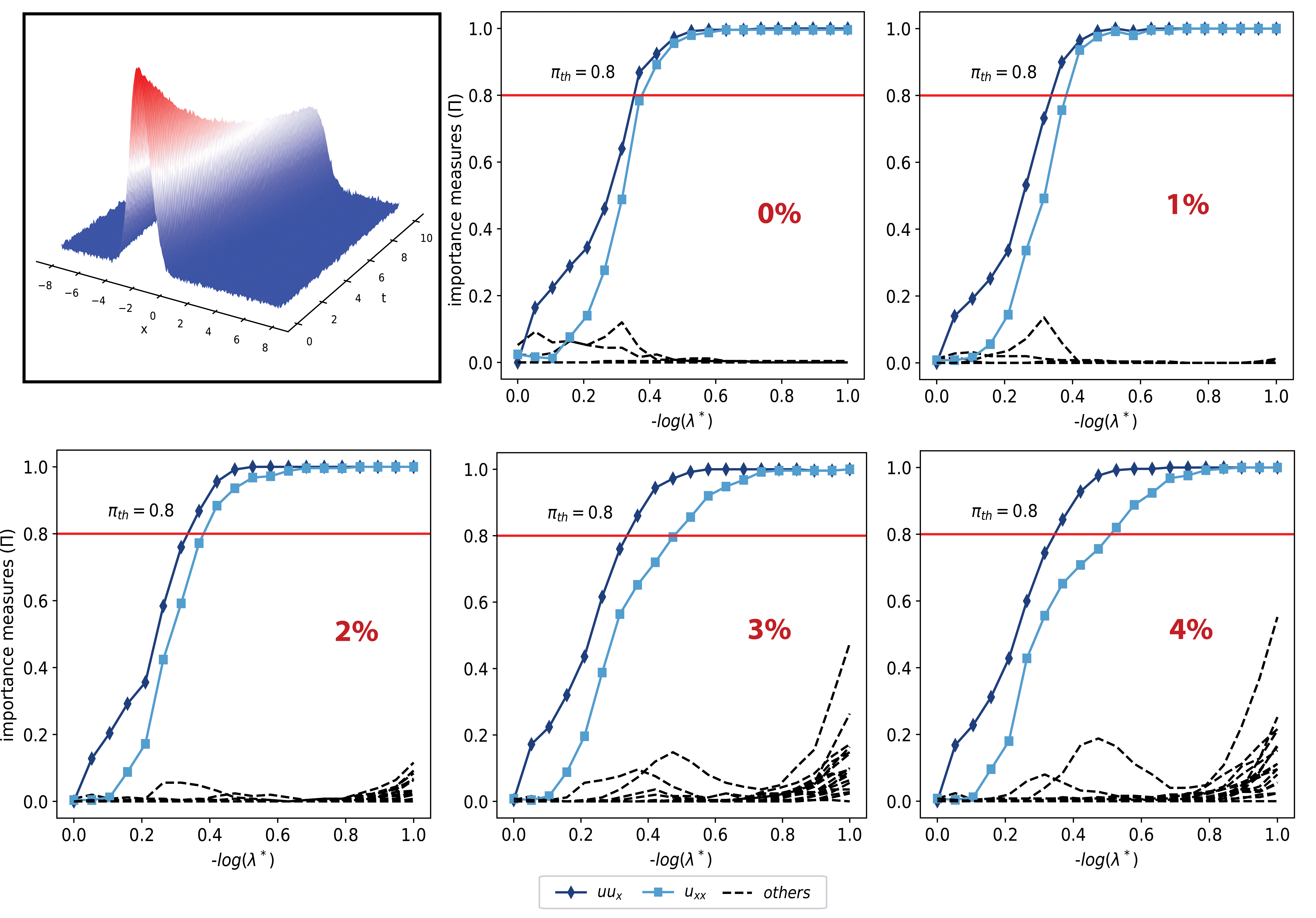}
\end{minipage}
\captionsetup{justification=raggedright,margin=0.5cm}
\caption{ \textbf{Model selection with PDE-STRIDE+IHT-d for 1D Burgers equation recovery :}   The top left image shows the numerical solution of the 1D Burgers equations on $256 \times 100$ space and time grid. The stability plots for the design $N=20,p=19$ show the separation of the true PDE components (in solid color) from the noisy components (dotted black). The inference power of the PDE-STRIDE method is tested for additive Gaussian noise-levels $\sigma$ up-to $5\%$ (not shown). In all the cases, perfect recovery was possible with the fixed threshold of $\pi^{th} = 0.8$ on the importance measure $\Pi$ (shown by the horizontal red solid line). The inset at the bottom shows the colors correspondence with the dictionary components.}
\label{Burgers_results_p19}
\end{figure}

We again consider the 1D Burgers equation from Eq.~(\ref{Burgers}), using the same simulated data as in Section \ref{result_1}, to quantify the performance and robustness against noise of the PDE-STRIDE+IHT-d method. The results are shown in Figure \ref{Burgers_results_p19} for a design with $N=250$ and $p = 19$. Even on this small data set, with a sample size comparable to dictionary size, our method recovers the correct model ($\{ u_{xx}, u u_{x} \}$) with up to $5 \%$ noise on the data, although the least-squares fits of the coefficient values gradually deviate from their exact values (see Table \ref{coefficient_values_Burgers}).

\noindent For comparison, the corresponding stability plots for PDE-STRIDE+STRidge are shown in Supplementary Figure~\ref{Burgers_results_STR_p20}. When using STRidge regression, the algorithm creates many false positives, even at mild noise levels ($<2$\%).

\begin{table}[h]
    \centering
    \scalebox{0.9}{
\begin{tabular}{ |c|c|c| } 
 \hline
  & $u u_{x} (-1.0) $ & $ u_{xx} ( 0.1 )$ \\
 \hline
 clean & -1.0008 & 0.1000 \\ 
 \hline
 $1 \% $ & -0.9971 & 0.1016 \\ 
 \hline
 $2 \% $ & -0.9932 & 0.0997 \\
 \hline
  $3 \% $ & -0.9842 & 0.0976 \\ 
 \hline
  $4 \% $ & -0.9728 & 0.0984 \\ 
 \hline
  $5 \% $ & -0.9619 & 0.0967 \\ 
 \hline
\end{tabular}}
\vspace{0.5em}
\captionsetup{justification=raggedright,margin=1.0cm}
\caption{  \footnotesize  Coefficients values of the recovered 1D Burgers equation for different noise levels.The stable components of the PDE inferred from plots in Figure \ref{Burgers_results_p19} are $\hat{S}_{stable} = \{ u_{xx}, uu_{x} \}$. }
\label{coefficient_values_Burgers}
\end{table}

\subsubsection{2D vorticity transport equation}
\begin{figure}[h]
\begin{minipage}[b]{1.0\textwidth}
\centering
\includegraphics[width=0.75\linewidth]{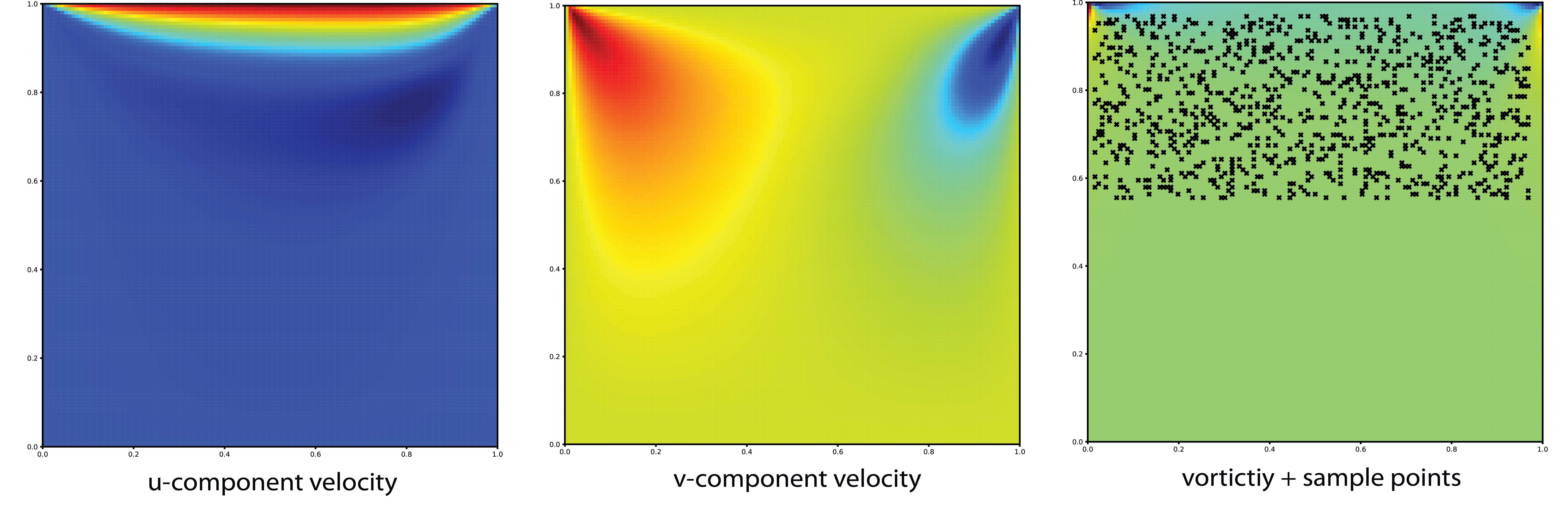}
\end{minipage}
\captionsetup{justification=raggedright,margin=0.5cm}
\caption{ \textbf{Numerical solution of 2D Vorticity transport equation}: The 2D domain with $u,v$ velocity components and vorticity $\omega$ is illustrated in the square domain. We choose to sample inside the rectangular box $ [0,1.0] \times [0.6,1.0] $ in the upper part of the domain to capture the rich dynamics resulting from shear boundary conditions imposed at the top surface. The black dots ($\approx 2000$) denote the points at which the data is sampled.}
\label{NS_illustrate}
\end{figure}

\begin{figure}[h]
\begin{minipage}[b]{1.0\textwidth}
\centering
\includegraphics[width=0.75\linewidth]{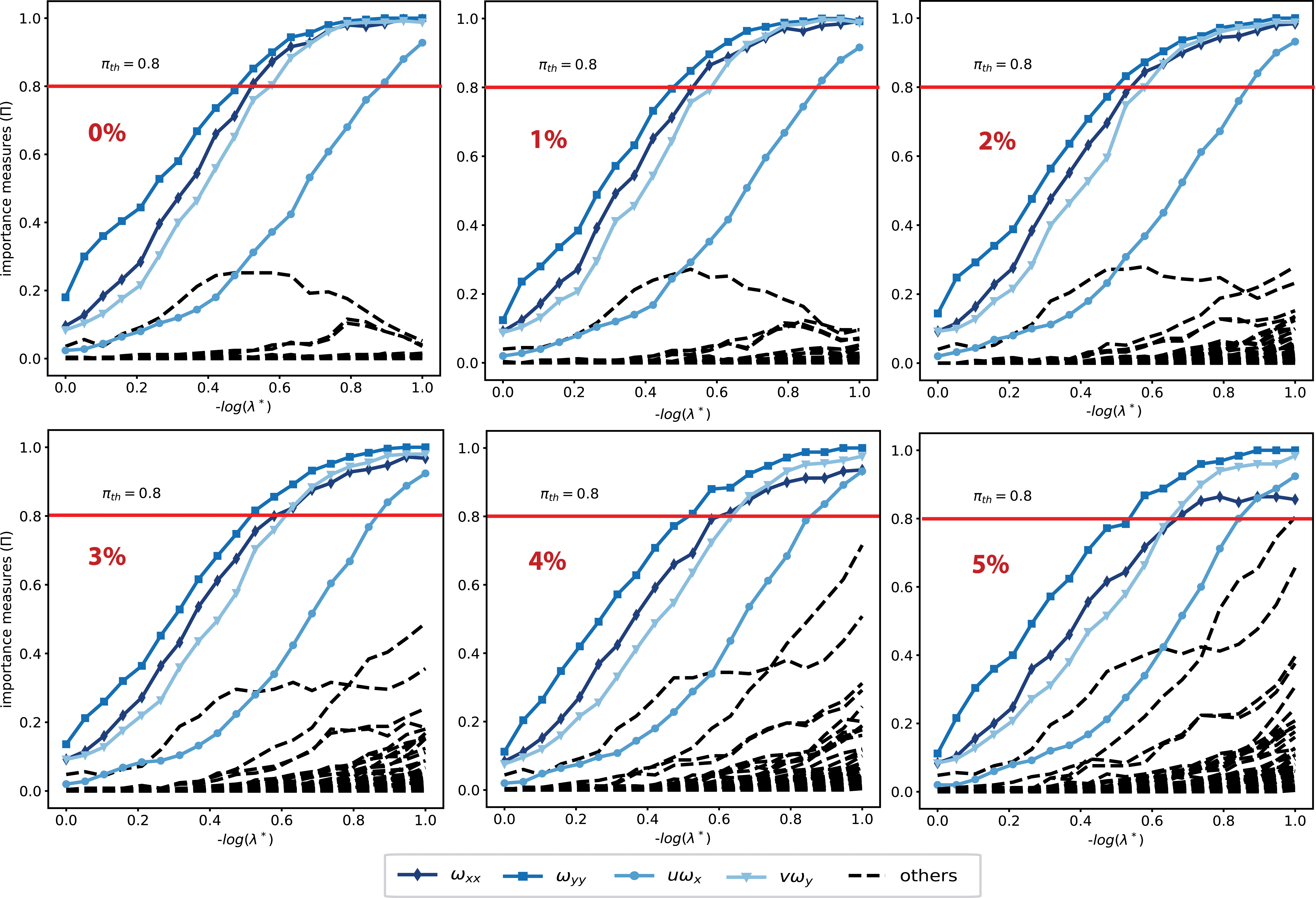}
\end{minipage}
\captionsetup{justification=raggedright,margin=0.5cm}
\caption{ \textbf{Model selection with PDE-STRIDE+IHT-d for 2D Vorticity transport equation recovery:} The stability plots for the design $N=500, p=48$ show the separation of the true PDE components (in solid color) from the noisy components (dotted black). The inference power of the PDE-STRIDE method is tested for additive Gaussian noise-levels $\sigma$ up-to $6\%$ (not shown). In all the cases, perfect recovery was possible with the fixed threshold of $\pi^{th} = 0.8$ on the importance measure $\Pi$ (shown by the horizontal red solid line). The inset at the bottom shows the colors correspondence with the dictionary components.}
\label{NS_results_p48}
\end{figure}

This section discusses results for the recovery of 2D vorticity transport equation using PDE-STRIDE. The vorticity transport equation can be obtained by taking curl of the Navier-Stokes equations and imposing the divergence-free constraint for enforcing in-compressibility, i.e. $\nabla \cdot u = 0$. This form of Navier-Stokes has found extensive applications in oceanography and climate modeling \cite{temam2001navier}. For the numerical solution of the transport equation, we impose no-slip boundary condition at the left $(x = 0, y \in [0,1])$, right $(x = 1, y \in [0,1])$ and bottom sides $(y = 0, x \in [0,1])$ and shear flow boundary condition $U = 2.0, V = 0$ on the top side $(y = 1, x \in [0,1])$. The simulation code was written using openFPM framework \cite{incardona2019openfpm} with explicit-time stepping on a $128 \times 128$ space grid. The poisson problem was solved at every time-step to correct velocities $u,v$ to ensure divergence-free fields. The viscosity of the fluid simulated was set to $\mu = 0.025$. In Figure \ref{NS_illustrate}, we show a single time snapshot of the $u,v,$ velocities and the vorticity field $\omega$ inside the square domain $[0,1] \times [0,1]$ of the Lid-driven cavity experiment.
\vspace{-0.5em}
\begin{equation} \label{NS}
    \omega_t + u \omega_x + v \omega_y = \mu \left( \omega_{xx} + \omega_{yy} \right)
\end{equation}
In Figure \ref{NS_results_p48}, the PDE-STRIDE results for 2D vorticity transport equation recovery are shown. The results demonstrate consistent recovery of the true support of the PDE for different noise-levels $\sigma$. The stable components $\hat{S}_{stable} = \{ \omega_{xx}, \omega_{yy}, u\omega_{x}, v\omega_{y} \}$ recovered correspond to the true PDE components of the 2D vorticity equation Eq. \ref{NS}. In table \ref{coefficient_values_NS}, we show refitted coefficients for the recovered PDE components. It should also be noted that the separation between the true (colored solid-lines) and the noisy (black dotted-lines) becomes less distinctive with increasing noise-levels. In the supplementary Figure \ref{NS_STR_results_p48}, we also report the STRidge based stability selection results for the same design and stability selection parameters. It can be seen that STRidge struggles to recover the true support even at small noise-levels, i.e. $\sigma > 0.01$.

\begin{table}[h]
\centering
\scalebox{0.9}{
\begin{tabular}{ |c|c|c|c|c| } 
 \hline
  & $ \omega_{xx} \: (0.025) $ & $\omega_{yy} \: (0.025) $ & $u\omega_{x} \: (-1.0) $ & $ v \omega_{y} (-1.0) $ \\
 \hline
 clean & 0.02504 & 0.02502 & -0.9994 & -1.0025 \\ 
 \hline
 $1 \% $ & 0.02501 & 0.02504 & -0.9997 & -1.0006\\ 
 \hline
 $2 \% $ & 0.02492 & 0.0250 & -1.0003 & -0.9944 \\
 \hline
  $3 \% $ & 0.0247 & 0.0250 & -1.004 & -0.9841 \\ 
 \hline
  $4 \% $ & 0.0245 & 0.0251 & -1.0091 & -0.9748 \\ 
 \hline
  $5 \% $ & 0.0242 & 0.0251 & -1.0083 & -0.9586\\ 
 \hline
\end{tabular}}
\vspace{0.5em}
\captionsetup{justification=raggedright,margin=1.0cm}
\caption{ \footnotesize Coefficients of the recovered 2D Vorticity transport equation for different noise levels. The stable components of the PDE from Figure \ref{NS_results_p48} are $\hat{S}_{stable} = \{ \omega_{xx}, \omega_{yy}, u\omega_{x}, v \omega_{y} \} $. }
\label{coefficient_values_NS}
\end{table}

\subsubsection{3D Gray-Scott equation}
In this section, we report the recovery capabilities of PDE-STRIDE for the 3D Gray-Scott reaction-diffusion equation Eqs. \ref{subequation_GS}. Reaction and diffusion of chemical species can produce a variety of patterns, reminiscent of those often observed in nature. Such processes form essential basis for morphogenesis in biology \cite{meinhardt1982models} and may even be used to describe skin patterning and pigmentation \cite{manukyan2017living}. We choose this example to show examples of systems with coupled variables and is very similar in dynamics to our real world example discussed in section \ref{result_5}. The reaction-diffusion dynamics described by Eqs. \ref{subequation_GS} are commonly be used to model the non-linear interactions between two chemical species $(u, v)$.
\begin{subequations}\label{subequation_GS}
\begin{align}
 u_t &= D_u \left ( \frac{\partial^2 u}{\partial x^2} +\frac{\partial^2 u}{\partial y^2} + \frac{\partial^2 u}{\partial z^2} \right) - uv^2 + f(1-u),\\
    v_t &= D_v \left( \frac{\partial^2 u}{\partial x^2} +
    \frac{\partial^2 u}{\partial y^2} + \frac{\partial^2 u}{\partial z^2} \right) + uv^2 - (f+k)v.
\end{align}
\end{subequations}

\begin{figure}[h]
\begin{minipage}[b]{1.0\textwidth}
\centering
\includegraphics[width=0.75\linewidth]{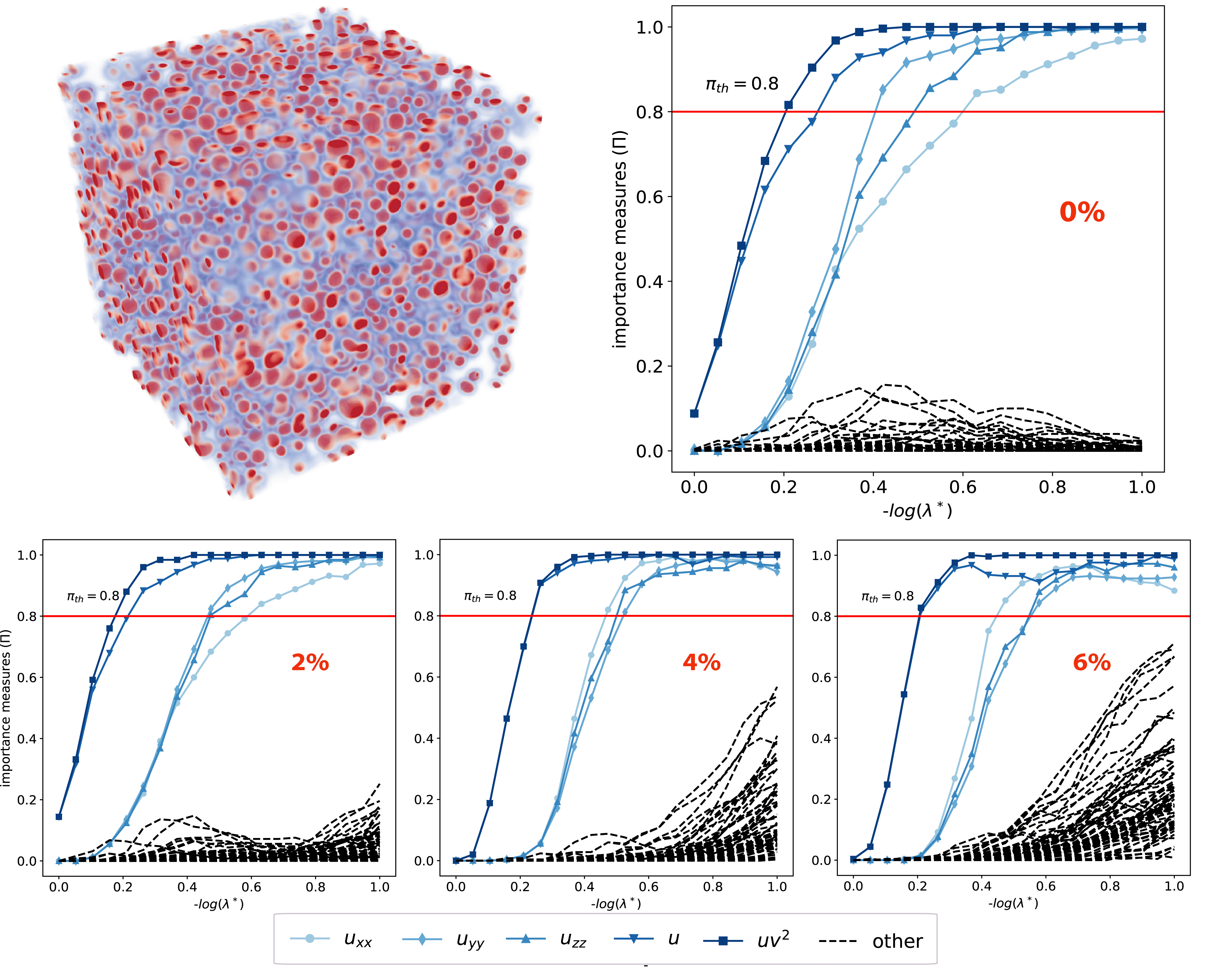}
\end{minipage}
\captionsetup{justification=raggedright,margin=0.5cm}
\caption{ \textbf{Model selection with PDE-STRIDE+IHT-d for 3D Gray-Scott $u-$component
equation recovery :}  The top left figure shows the visualization of the 3D simulation domain with $v$ species concentration. The color gradient corresponds to the varying concentration over space. The stability plots for the design $N=400, p=69$ show the separation of the true PDE components (solid color) from the noisy components (dotted black). The inference power of the PDE-STRIDE method is tested for additive Gaussian noise-levels $\sigma$ up-to $6\%$ (not shown). In all the cases, perfect recovery was possible with the fixed threshold of $\pi^{th} = 0.8$ on the importance measure $\Pi$ (shown by the horizontal red solid line). The inset at the bottom shows the colors correspondence with the dictionary components.}
\label{Fig_gray_scott_u}
\end{figure}

\noindent We simulate the above equations using openFPM framework. The snapshot of the 3D cube simulation box $2.5 \times 2.5 \times 2.5$ along with the $v$ concentration is shown in Fig. \ref{Fig_gray_scott_u}. Finite-difference space discretization scheme was used to discretize the dynamics described in Eqs. \ref{subequation_GS} with grid spacing $ \textrm{dx} = \textrm{dy} = \textrm{dz} = 0.01953 $. The explicit time-stepping with step size $\textrm{dt} = 0.0005$ was used for temporal integration to simulate for $5 s$ in real time. The 3D Gray-Scott model parameters used are $k = 0.053, f = 0.014, D_{u} = 2.0\textrm{E}^{-5} $, and $D_{v} = 1.0\textrm{E}^{-5}$. \\

\noindent Given the large degrees of freedom present in the 3D problem, for our learning problem we choose to sample data only from a small cube in the middle of the domain with dimension $0.5 \times 0.5 \times 0.5$. In Figure \ref{Fig_gray_scott_u}, we show PDE-STRIDE method correctly identifies the true PDE components for the dynamics of $u$ species given by Eq. (3.3.2a) for different noise levels. One can appreciate the clear separation between the true and noisy PDE components in the stability plots. We show results for different noise-levels between $0-6 \%$ with as few as $N = 400$ samples and for dictionary size $p = 69$. Similar plots for the inference of $v$ species dynamics are shown in Fig. \ref{Fig_gray_scott_v}.  Although perfect recovery was not possible owing to the small diffusivity ($D_v = 1.0\textrm{E}^{-5}$) of the $v$ species, consistent and stable recovery of the reaction terms (interaction terms) can be seen. The refitted coefficients for the recovered PDE for the both the $u,v$ species are reported in table \ref{coefficient_values_uGS} and table \ref{coefficient_values_vGS}, respectively. The comparison plots for PDE-STRIDE with STRidge for the 3D Gray-Scott recovery are shown in the supplementary Figures \ref{uGS_STR_results_p48}, \ref{vGS_STR_results_p48}. We note that the STRidge is able to recover the complete form of Eq \ref{subequation_GS} in noise-free case for both the $u,v$ species, but it fails to recover both the $u,v$ PDEs in the noise case. The comparison clearly demonstrates that PDE-STRIDE+IHT-d  clearly outperforms PDE-STRIDE+STRidge for inference from noisy data-sets.

\begin{table}[h]
\centering
\scalebox{0.9}{
\begin{tabular}{ |c|c|c|c|c|c|c| } 
 \hline
  & 1(0.014) &$ u_{xx} (2.0\textrm{E}^{-5}) $ & $u_{yy} (2.0\textrm{E}^{-5}) $ & $u_{zz} (2.0\textrm{E}^{-5}) $ & $ u (-0.014) $ & $ uv^2 (-1.0) $ \\
 \hline
 clean & 0.0140 & $2.0\textrm{E}^{-5}$ & $2.0\textrm{E}^{-5}$ & $2.0\textrm{E}^{-5}$ & -0.0140  & -1.0000\\ 
 \hline
 $2 \% $ & 0.0142 & $1.9664\textrm{E}^{-5}$ & $1.9565\textrm{E}^{-5}$ & $1.9869\textrm{E}^{-5}$ & -0.0143 & -0.9915 \\ 
 \hline
 $4 \% $ & 0.0144 & $1.9541\textrm{E}^{-5}$ & $1.8971\textrm{E}^{-5}$ & $1.8780\textrm{E}^{-5}$ & -0.0146 & -0.9795\\
 \hline
  $6 \% $ & 0.0150 & $2.0494\textrm{E}^{-5}$ & $1.8888\textrm{E}^{-5}$ & $1.8284\textrm{E}^{-5}$ &  -0.0153 &-0.9843\\ 
 \hline
\end{tabular}}
\vspace{0.5em}
\captionsetup{justification=raggedright,margin=1.0cm}
\caption{ Coefficients of the recovered $u$-component Gray-Scott reaction diffusion equation for different noise levels. The stable components of the PDE from Figure \ref{Fig_gray_scott_u} are $\hat{S}_{stable} = \{ u_{xx}, u_{yy}, u_{zz}, u,  uv^2 \}$  }
\label{coefficient_values_uGS}
\end{table}

\noindent In all the above presented cases, the learned stable components $(\hat{S}_{stable})$ coincide with the true components of the underlying PDE model. The PDE-STRIDE framework is able to learn the stable components with data points as few as $\approx 400$. So, we conclude that our PDE-STRIDE+IHT-d framework is able to robustly learn partial differential equations from limited noisy data. In the next section, we discuss the consistency and robustness of the PDE-STRIDE+IHT-d for perturbations in design parameters like sample-size $N$, dictionary size $p$, and noise-levels $\sigma$.

\subsection{Achievability results} \label{result_4}
Achievability results are a compact way to check for the robustness and consistency of a model selection method for varying design parameters. They also provide an approximate means to reveal the \textit{sample complexity} of any $l_{0}$ and $l_{1}$ sparsity-promoting technique, i.e. the number of data points $N$ required to recover the model with full probability. Specifically, given a sparsity promoting method, dictionary size $p$, sparsity $k$, and noise level $\sigma$, we are interested in how the sample size $N$ scales with $p, k, \sigma$ with recovery probability converging to one. The study \cite{wainwright2009sharp} reports sharp phase transition from failure to success for Gaussian random designs with increasing sample size $N$ for LASSO based sparsity solutions. The study also provides sufficient lower bounds for sample size $N$ as a function of $p,k$ for full recovery probability. In this section, we ask the question on whether sparse model selection with PDE-STRIDE (paired with IHT-d) exhibit similar sharp phase transition behaviour. Given the dictionary components in our case are compiled from derivatives and non-linearities computed from noisy-data, it is interesting to observe if full recovery probability is achieved and maintained with increasing sample size($N$). In the particular context of PDE learning, increasing dictionary size by including higher order non-linearities and higher order derivatives has the potential to include strongly correlated components, which can negatively impact the inference power. This observation was made evident with results and discussion from section \ref{result_2}.  \\

\begin{figure}[h]
\begin{minipage}[b]{1.0\textwidth}
\centering
\includegraphics[width=0.85\linewidth]{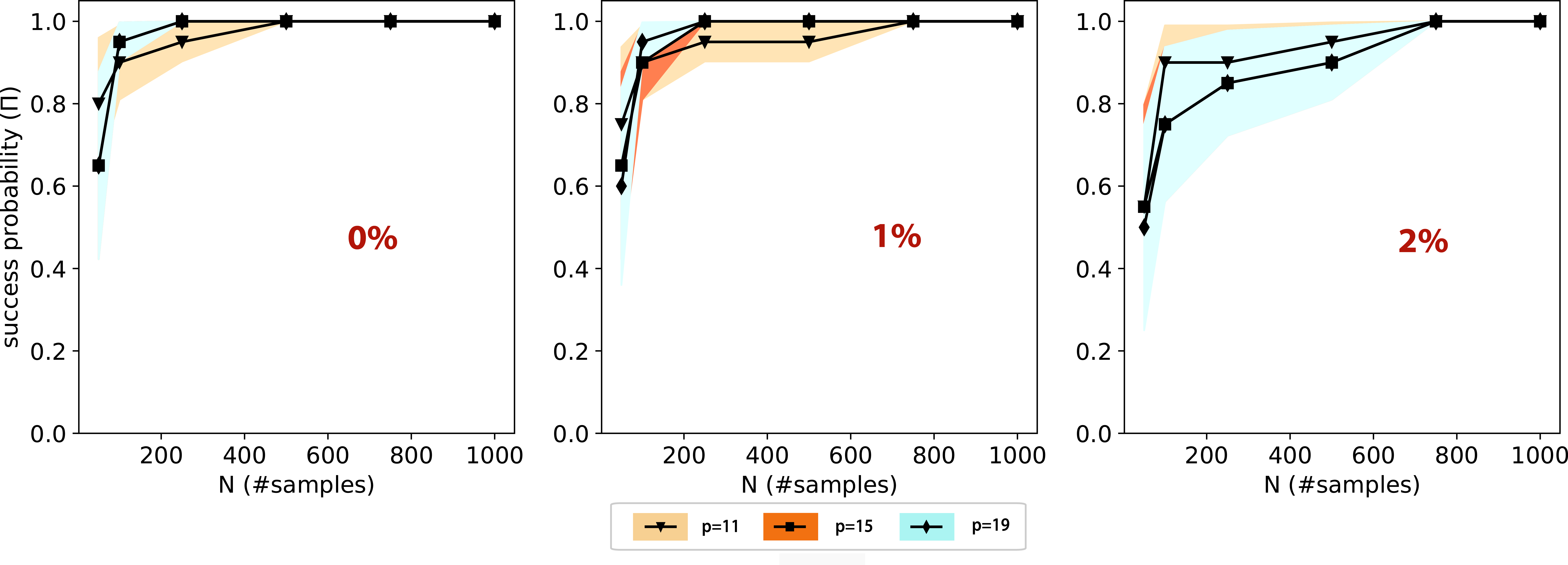}
\end{minipage}
\captionsetup{justification=raggedright,margin=0.5cm}
\caption{ \textbf{Achievability results for model selection with PDE-STRIDE+IHT-d for 1D Burgers equation recovery :} The achievability results for 1D Burgers equation recovery with PDE-STRIDE is shown for varying designs $(N, p, \sigma)$. Every point on the plot corresponds to 20 repetitions of the PDE-STRIDE method for some design $(N, p, \sigma)$ under random data sub-sampling. Each line with markers corresponds to a dictionary size $p$. In Burgers case, we test for $p = \{11, 15, 19\}$ clearly shown in the inset below the plots. The colored area around the line show the associated variance of the Bernoulli's trials. }
\label{Burgers_AP}
\end{figure}

\begin{figure}[h]
\begin{minipage}[b]{1.0\textwidth}
\centering
\includegraphics[width=0.85\linewidth]{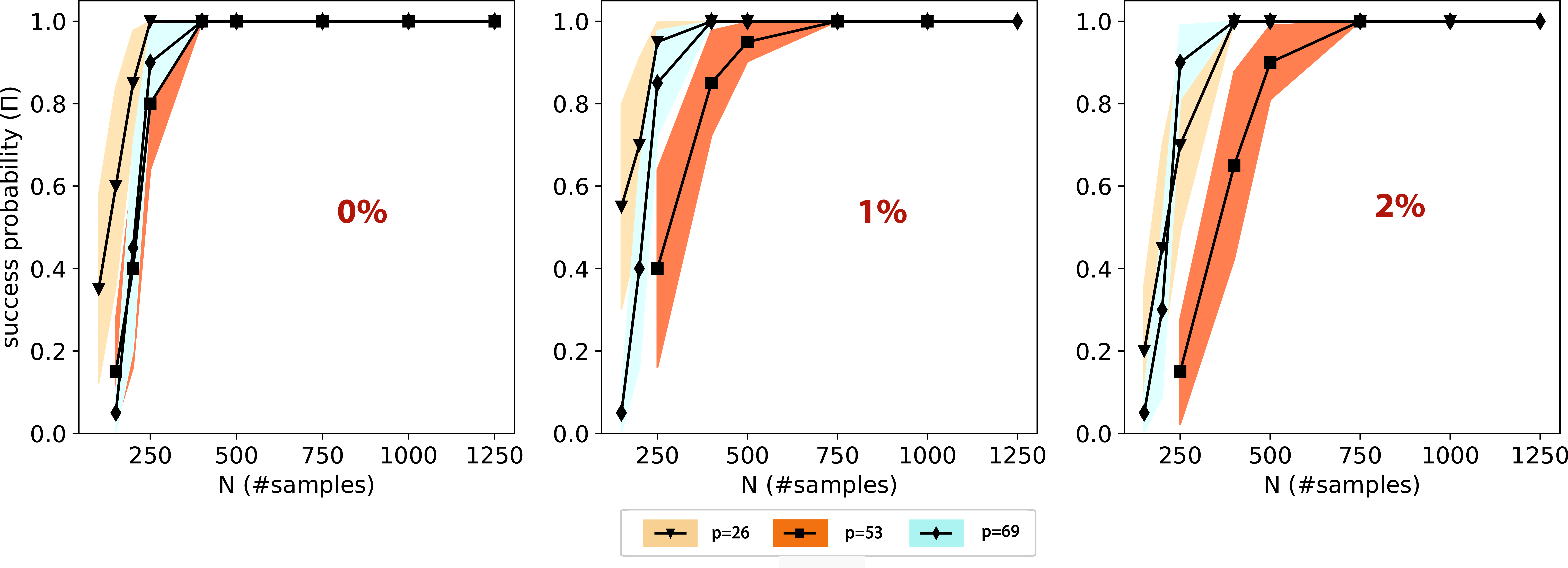}
\end{minipage}
\captionsetup{justification=raggedright,margin=0.5cm}
 \caption{ \textbf{Achievability results for model selection with PDE-STRIDE+IHT-d for 3D Gray-Scott $u$-component equation recovery :} The achievability results for 3D Gray-Scott equation recovery with PDE-STRIDE method is shown for varying designs $(N, p, \sigma)$. Every point on the plot corresponds to 20 repetitions of the PDE-STRIDE for some design $(N, p, \sigma)$ under random data sub-sampling. Each line with markers corresponds to a dictionary size $p$. In 3D Gray-Scott case, we test for $p = \{26, 53, 69\}$ clearly shown in the inset below the plots. The colored area around the line show the associated variance of the Bernoulli's trials. }
\label{GS_AP}
\end{figure}

\noindent In Figures (\ref{Burgers_AP}, \ref{GS_AP}) the achievability results for both the 1D Burgers system and 3D $u$-component Gray-Scott reaction diffusion system are shown. Each point in Figure (\ref{Burgers_AP}, \ref{GS_AP}) correspond to 20 repetitions of an experiment with some design $(N, p, \sigma)$ under random data sub-sampling. An experiment with a design $ (N, p, \sigma) $ is considered as success if $\exists	\lambda \in \Lambda$ for which the true PDE support $(S^{*})$ is recovered with PDE-STRIDE+IHT-d by thresholding at $\pi_{th}=0.8$. In Figures (\ref{Burgers_AP}, \ref{GS_AP}), we see strong consistency and robustness to design parameters for both Burgers and Gray-Scott systems. And, we also observe sharp phase transition from failure to success with recovery probability converging to one with increasing sample size $N$. This amply evidence suggest that PDE-STRIDE not only enhances the inference power of the IHT-d method but also ensures consistency. In addition, the sharp phase transition behaviour also point towards a strict lower bound on the \textit{sample complexity} (N) below which full recoverability is not attainable as studied in \cite{wainwright2009sharp}. Following this, achievability plots can be used to read-out approximate estimates of the \textit{sample-complexity} of the learned dynamical systems. In the case of Burgers, $90 \%$ success probability is achieved with data points as few as $ \approx 70$ in noise-free and $\approx 200$ points in noisy cases for different designs ($p$). For the 3D Gray-scott system,  $90 \%$ success probability is achieved with data points as few as $ \approx 200$ in noise-free and $\approx 400$ points in noisy cases for different designs ($p$).
\vspace{-1.0em}

\begin{figure}[!h]
\begin{minipage}[b]{1.0\textwidth}
\centering
\includegraphics[width=0.85\linewidth]{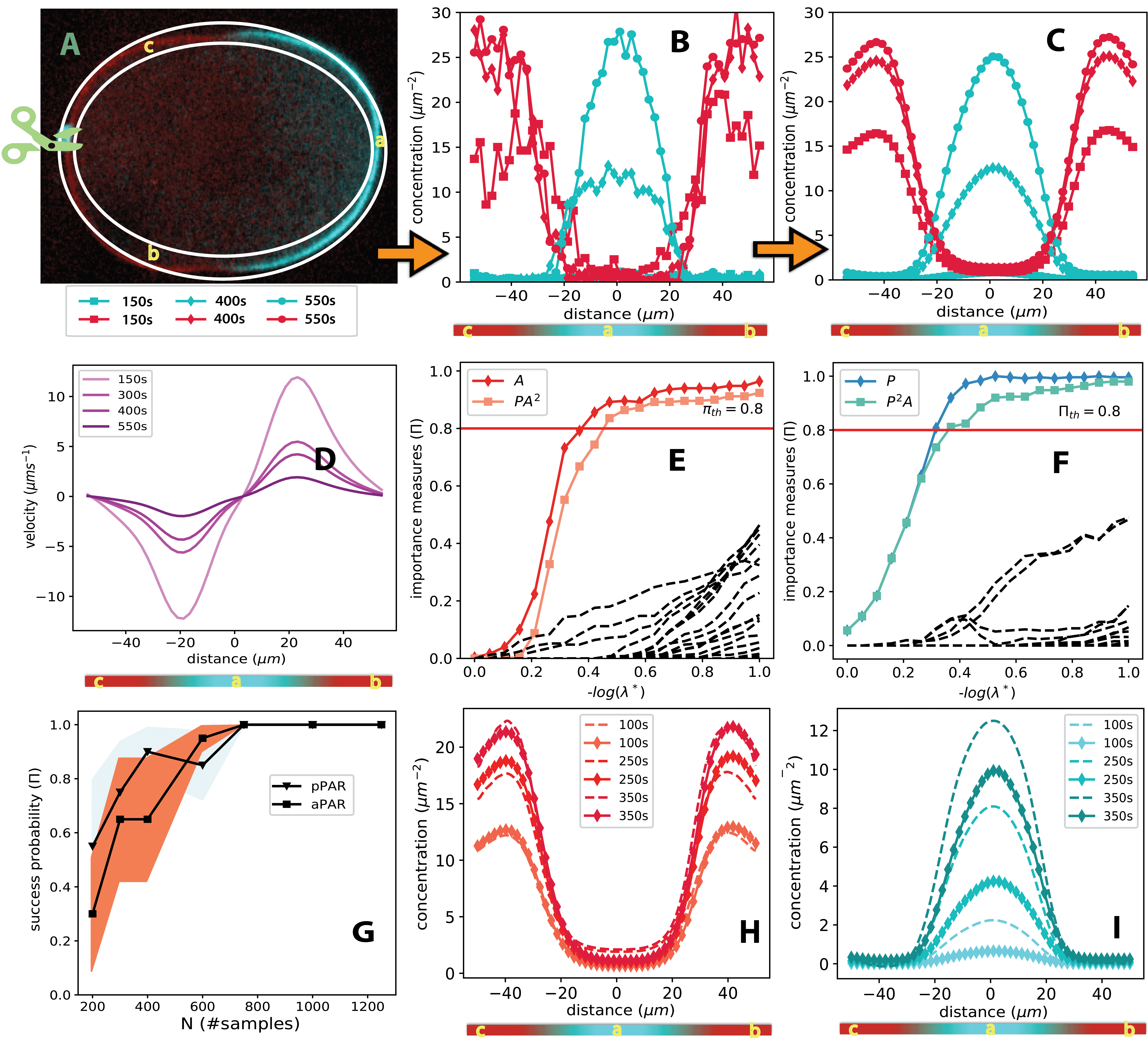}
\end{minipage}
\vspace{0.15em}
\captionsetup{justification=raggedright,margin=0.5cm}
\caption{ \textbf{ Data-driven model inference of the regulatory network of membrane PAR proteins from spatiotemporal data acquired from \textit{C. elegans} zygote :} \textbf{A}:  spatiotemporal PAR concentration data-sets provided by Grill Lab \cite{etemad1995asymmetrically,goehring2011polarization} at MPI-CBG. Manually defining observational boundaries (white ellipses) and identification of key variables of interest like protein concentration and cell-cortex velocity in the microscopy data. The blue color and red color corresponds to the intensities of the posterior PAR (pPAR) and anterior PAR (aPAR) proteins, respectively. \textbf{B}: The noisy aPAR(red) and pPAR(blue) concentration profiles extracted from the experiments from fluorescence microscopy. \textbf{C}: De-noised spatiotemporal concentration profiles obtained from extracting the principle mode of the Singular value decomposition (SVD) of the noisy flow data. Different symbol lines correspond to different time instances shown in the bottom inset of Figure A. \textbf{D}: De-noised spatiotemporal cortical flow profiles obtained from extracting the principle mode of the Singular value decomposition (SVD) of the noisy data. \textbf{E}: The stability plots using PDE-STRIDE+IHT-d to identify the stable PDE components (colored) and learn the model for aPAR protein interaction. \textbf{F}: The stability plots using PDE-STRIDE+IHT-d to identify the stable PDE components (colored) and learn the model for pPAR protein interaction. \textbf{G}: Achievability results to test the robustness and consistency of the inferred model of both aPAR - $ S_{stable}^{A} = \{ A, A^2 P \} $ and pPAR - $ S_{stable}^{P} = \{ P, P^2 A \} $ with increasing sample-size $N$. \textbf{H}: The simulation results ($\color{red}\dashedrule$) of the learned models for aPAR protein overlapped with the experimental data ($\color{red}-\!\! \vardiamondsuit \!\!-$) at different times.  \textbf{I}: The simulation results ($\color{cyan}\dashedrule$) of the learned model for pPAR protein overlapped with the experimental data ($\color{cyan}-\!\! \vardiamondsuit \!\!-$)  at different times.  }
\label{PAR}
\end{figure}

\section{Data-driven PDE inference on real experimental data to explain \textit{C.~elegans} zygote patterning} \label{result_5}

We showcase the applicability of the PDE-STRIDE with IHT-d to real experimental data. We use microscopy images to infer a PDE model that explains early {\it C.~elegans} embryo patterning, and we use it to confirm a previous hypothesis about the physics driving this biological process. Earlier studies of this process proposed a mechano-chemical mechanism for PAR protein polarization on the cell membrane~\cite{etemad1995asymmetrically,goehring2011polarization,gross2019guiding}. They systematically showed that the cortical flows provide sufficient perturbations to trigger polarization \cite{goehring2011polarization}.  The experiments conducted in \cite{goehring2011polarization} measured the concentration of the anterior PAR complex (aPAR), the concentration of posterior PAR complex (pPAR) and the cortical flow field ($\textrm{v}$) as a function of time as shown in Figure \ref{PAR}, A, B, D. The concentration and velocity fields were acquired on a grid with resolution of $60 \times 55 $ in space and time. These experimental data-sets were used to validate the mechano-chemical model developed in the studies \cite{etemad1995asymmetrically,goehring2011polarization}.  Here, we challenge the PDE-STRIDE+IHT-d framework to learn a PDE model for the regulatory network (Eq \ref{reaction_format}) of the interacting membrane PAR proteins in a pure data-driven sense from the experimental data-set. Given the noisy nature of the data-sets, our analysis is limited to the first SVD mode of the data as shown in Figure \ref{PAR} C. We also focus our attention on the temporal regime post the advection trigger when the early domains of PAR proteins are already formed. The PDE-STRIDE is then directed to learn an \textit{interpretable} model from the data, that evolves the early protein domains to fully developed patterns as shown in Figure \ref{PAR}A.\\

\noindent The reaction kinetics of the PAR proteins can be formulated as,
\vspace{-0.5em}
\begin{equation}\label{reaction_format}
    v_{a}^{-} A +  v_{p}^{-} P \xrightarrow{k} v_{a}^{+} A +  v_{p}^{+} P.
\end{equation}
Here, $v_{a/p}^{-}$ and $v_{a/p}^{+}$ are the reactant and product stoichiometry, respectively. The variables $A$ and $P$ correspond to the aPAR and pPAR protein species and $k$ is the reaction rate.

\noindent  In designing the dictionary $\Theta$, the maximum allowed stoichiometry for reactant and product is restricted to 2, i.e. $v_{a/p}^{-}, v_{a/p}^{+} \in \{ 0,1,2 \}$. The PDE-STRIDE+IHT-d results for the learned regulatory reaction network from data are shown in Figure \ref{PAR} E, F. The stable components of the model for aPAR protein are $\hat{S}_{stable}^{P} = \{ P, P^2 A \}$ and for pPAR protein are $ \hat{S}_{stable}^{A} = \{ A, PA^2\}$ for a design $N\approx500, p=20$. In Figure \ref{PAR}G, achievability tests are conducted to show the consistency and robustness of the learned models across different sample-sizes $N$. The learned model achieves full recovery probability for sample-size $N > 800$. Our preliminary models inferred in a data-driven manner exhibit very good qualitative agreement with the experiments and also recapitulates the mutual inhibitory nature of the PAR protein. The parameters of the learned models are then computed by least-squares refitting and are tabulated in Table \ref{coefficient_values_PAR}.
\begin{table}[h]
    \centering
    \scalebox{0.8}{
\begin{tabular}{ |c|c|c|c| } 
 \hline
pPAR & $1 $ &$ P $ & $ AP^2$ \\
 \hline
 &-0.00019769 & 0.01073594 & -0.00027887\\ 
 \hline
aPAR & $1 $ &$ A $ & $ A^2P$ \\
 \hline
 &  0.0041325 & -0.00216077 & -0.00014699\\ 
 \hline
\end{tabular}}
\vspace{0.5em}
\caption{  \footnotesize  Coefficients values of the inferred PAR model. The stable components of the PDE inferred from stability results in Figure \ref{PAR} E, \ref{PAR} F  are $\hat{S}_{stable}^{P} = \{ P, P^2A \}$ and $\hat{S}_{stable}^{A} = \{ A, A^2P \}$. }
\label{coefficient_values_PAR}
\end{table}

\noindent In the Figure \ref{PAR}H, \ref{PAR}I, we overlay the numerical solution of the the learned model with the de-noised experimental data at certain time snapshots for both models of aPAR and pPAR proteins for quantitative comparison. This very simple PDE model is able to describe the temporal evolution of the PAR protein domains on the \textit{C. elegans} zygote. Although, there is a good match in the spatial-scales (PAR domain sizes) for the two proteins, there exists a non-negligible discrepancy between the simulation and experiments in the time-scales for the pPAR evolution. This difference can be attributed to the advection processes dictated by the cortical flows, which are not included in our simple ODE type model as shown in \ref{PAR}E, \ref{PAR}I. We believe including higher modes of the SVD decomposition and also using structured sparsity for enforcing symmetric arguments through grouping \cite{ward1963hierarchical} can further mature our data-driven models to include the mechanical aspects of the PAR system. In supplementary Figure \ref{extra_PAR}(left), we already show for a particular design of $N=500, p=20$, the advection and diffusion components of the aPAR protein model carry significant importance measure to be included in the stable set $\hat{S}_{stable}$, but this is not the case with the pPAR components as shown in Figure \ref{extra_PAR}(right). The preferential advective displacement of the aPARs to the anterior side modeled by the advective term $(\textrm{v}_x A)$ is also in line with the observations of the experimental studies \cite{goehring2011polarization}. However, such models with advection and diffusion components exhibit inconsistency for varying sample-size $N$, in contrast to our simple ODE type model as illustrated in Figure \ref{PAR}G. 

\section{Conclusion and Discussion} \label{conclusion}
We have addressed two key issues that have thus far limited the application of sparse regression methods for automated PDE inference from noisy and limited data: the need for manual parameter tuning and the high sensitivity to noise in the data. We have shown that stability selection combined with any sparsity-promoting regression technique provides an appropriate level of regularization for consistent and robust recovery of the correct PDE model. Our numerical benchmarks suggested that iterative hard thresholding with de-biasing (IHT-d) is an ideal combination with  stability selection to form a robust and parameter-free framework (PDE-STRIDE) for PDE learning. This combination of methods outperformed all other tested algorithmic approaches with respect to identification performance, amount of data required, and robustness to noise. The resulting stability-based  PDE-STRIDE method was tested for robust recovery of the 1D Burgers equation, 2D vorticity transport equation, and 3D  Gray-Scott reaction-diffusion equations from simulation data corrupted with up to 6$\%$ of additive Gaussian noise. The achievability studies  demonstrated the consistency and robustness of the PDE-STRIDE method for full recovery probability of the model with increasing sample size $N$ and for varying dictionary size $p$ and noise levels $\sigma$. 
In addition, we note that achievability plots provide a natural estimate for the \textit{sample-complexity} of the underlying non-linear dynamical system. However, this empirical estimate of the \textit{sample-complexity} depends on the choice of model selection algorithm and how the data is sampled.\\

\noindent We demonstrated the capabilities of the PDE-STRIDE+IHT-d framework by applying it to learn a PDE model of embryo polarization directly from fluorescence microscopy images of \textit{C.~elegans} zygotes. The model recovered the regulatory reaction network of the involved proteins and their spatial transport dynamics in a pure data-driven manner, with no knowledge used about the underlying physics or symmetries. The thus learned, data-derived PDE model was able to correctly predict the spatiotemporal dynamics of the embryonic polarity system from the early spatial domains to the fully developed patterns as observed in the polarized \textit{C.~elegans} zygote. The model we inferred from image data using our method confirms both the structure and the mechanisms of popular physics-derived cell polarity models. Interestingly, the mutually inhibitory interactions between the involved protein species, which have previously been discovered by extensive biochemical experimentation, were automatically extracted from the data. \\

\noindent 
Besides rendering sparse inference methods more robust to noise and parameter-free, stability selection has the important conceptual benefit of also providing interpretable probabilistic importance measures for all model components. This enables modelers to construct their models with high fidelity, and to gain an intuition about correlations and sensitivities. Graphical inspection of stability paths provides additional freedom for user intervention in semi-automated model discovery from data.\\

\noindent We expect that statistical learning methods have the potential to enable robust, consistent, and reproducible discovery of predictive and interpretable models directly from  observational data. Our parameter-free PDE-STRIDE framework provides a first step toward this, but many open issues remain. First, numerically approximating time and space derivatives in the noisy data is a challenge for noise levels higher than a few percent. This limits the noise robustness of the overall methods, regardless of how robust the subsequent regression method is. The impact of noise becomes even more severe when exploring models with higher-order derivatives or stronger non-linearities. Future work should focus on formulations that are robust to the choice of different discretization methods, while providing the necessary freedom to impose structures on the coefficients. Second, a principled way to constrain the learning process by physical priors, such as  conservation laws and symmetries, is lacking. Exploiting such structural knowledge about the dynamical system is expected to greatly improve learning performance. It should therefore be explored whether, e.g., structured sparsity or grouping constraints \cite{kutz2017data, ward1963hierarchical} can be adopted for this purpose. Especially in coupled systems, like the  Gray-Scott reaction-diffusion system and the PAR-polarity model, known symmetries could be enforced through structured grouping constraints.
\\

\noindent 
In summary, we believe that data-driven model discovery has tremendous potential to provide novel insights into complex processes, in particular in biology. It provides an effective and complementary alternative to hypothesis-driven approaches. We hope that the stability-based model selection method PDE-STRIDE presented here is going to contribute to the further development and adoption of these approaches in the sciences.\\

\noindent \textbf{Code and data availability:} The github repository for the codes and data can be found at \url{https://github.com/SuryanarayanaMK/PDE-STRIDE.git}.

\subsection*{Acknowledgements}
This work was in parts supported by the Deutsche Forschungsgemeinschaft (DFG, German Research Foundation) under Germany's Excellence Strategy – EXC-2068 – 390729961. We are grateful to the Grill lab at MPI-CBG/TU Dresden for providing the  spatiotemporal PAR concentration and flow field data and allowing us to use them in our showcase. We thank Nathan Kutz (University of Washington) and his group for making their code and data public.

\bibliographystyle{unsrt}
\bibliography{stability}

\newpage
\section{Supplementary Material}

\subsection{Algorithm}
\begin{algorithm}
\begin{multicols}{2}
\begin{algorithmic}
\footnotesize
  \State \!\!\!\!\!\!\!\textrm{Problem:} $ \Large \hat{\xi} = \arg\min_{\xi} \Vert U_t - \Theta \xi \Vert_{2}^{2} + \lambda \Vert \xi \Vert_{0}$
\vspace{0.5em}
  \State \!\!\!\!\!\!\!\textbf{IHD-d}$(\Theta, U_t, \lambda, \textrm{maxit} , \textrm{subit} )$: \vspace{0.25em}
  \State \textbf{Initialize:} $\xi^{0}=0$ \vspace{0.25em}
    \For{ $ n  = 0 \textbf{ to } \textrm{maxit}$}
        \State $\nabla g = -\Theta^{T}(U_t - \Theta \xi^{n}) $ \vspace{0.25em}
        \State $u^{n+1,1} = H_{\lambda}\left ( \xi^{n} - \mu^{1} \cdot \nabla g  \right), \: S^{n+1} = \textrm{supp}(u^{n+1,1}) $ 
        \For{$ l = 1  \textbf{ to } \textrm{subit}$} \vspace{0.25em} 
            \State $\nabla g_{s}=  \left ( -\Theta^{T}(U_t - \Theta u^{n+1,l}) \right)_{S^{n+1}} $ \vspace{0.25em}
            \State $u^{n+1, l+1} = \left( u^{n+1,l} - \mu^2 \cdot \nabla g_{s} \right)_{S^{n+1}} $ \vspace{0.25em}
            \If { $ \left(  \:\: \Vert U_t - \Theta u^{n+1, l+1} \Vert_{2}^{2} \leq \lambda \vert S^{n+1} \vert \:\:   \right) $ } \vspace{0.25em}
                \State $\textrm{return} \:\: \hat{\xi} = u^{n+1,l+1}$ \vspace{0.25em}
            \EndIf \vspace{0.25em}
        \EndFor 
         \State $\xi^{n+1} = u^{n+1,l+1}$ \vspace{0.25em}
  \EndFor \vspace{0.25em}
  \State \textrm{return} $\hat{\xi} = u^{n+1,l+1}$

  \columnbreak
  
   \State \!\!\!\!\!\!\!\textbf{IHT}$(\Theta, U_t, K, \textrm{maxit} )$: \vspace{0.25em}
   \State \textbf{Initialize:} $\xi^{0}=0$ \vspace{0.25em}
       \For{ $ n  = 0 \textbf{ to } \textrm{maxit}$}
             \State $\nabla g = -\Theta^{T}(U_t - \Theta \xi^{n}) $ \vspace{0.25em}
            \State $u^{n+1,1} =  H_{\lambda}\left ( \xi^{n} - \mu^{1} \cdot \nabla g  \right) $  \vspace{0.25em}
       \EndFor \vspace{0.25em}\\
  \noindent \!\!\!\!\!\!\!  \rule{8cm}{0.4pt}
  
\State \!\!\!\!\!\!\!\textrm{Problem:}  $ \hat{\xi} = \arg\min_{\xi} \Vert U_t - \Theta \xi \Vert_{2}^{2},  \:\:\: \textrm{s.t.} \:\:\:  \Vert \xi \Vert_{0} \leq K$
\vspace{0.5em}
  \State \!\!\!\!\!\!\!\textbf{HTP}$(\Theta, U_t, K, \textrm{maxit} )$: \vspace{0.25em}
   \State \textbf{Initialize:} $\xi^{0}=0$ \vspace{0.25em}
       \For{ $ n  = 0 \textbf{ to } \textrm{maxit}$}
             \State $\nabla g = -\Theta^{T}(U_t - \Theta \xi^{n}) $ \vspace{0.25em}
            \State $u^{n+1,1} =  \textrm{indices of K largest entries of } \xi^{n} - \mu^{1} \cdot \nabla g  $ \vspace{0.25em}
            \State $ S^{n+1} = \textrm{supp}(u^{n+1,1}) $  \vspace{0.25em}
            \State $\xi^{n+1}  =  \arg \min \bigg(  \Vert U_t - \Theta z \Vert \bigg)_{S^{n+1}}  $  \vspace{0.25em}
       \EndFor \vspace{0.25em}

\end{algorithmic}
\caption{$\hat{\xi} = \arg\min_{\xi} \Vert U_t - \Theta \xi \Vert_{2}^{2} + \lambda \Vert \xi \Vert_{0}$}
\label{algorithm_1}
\end{multicols}
\end{algorithm}

\noindent In the above algorithm $\mu^{1}$ and $\mu^{2}$ correspond to the step-size/learning rates corresponding to the IHT and de-biasing steps respectively. The learning rate $\mu_{1}$ is computed as the inverse of the Lipschitz constant of the gradient of the square-loss function $h(\cdot)$ in Eq.(\ref{square_loss}), i.e. $\mu_{1} = 1.0/L$. In the similar fashion, the learning rate in the de-biasing step is computed as $\mu_{2} = 1.0/L^{*}$. Here, $L^{*}$ is the Lipschitz constant of the square-loss function $ \left (h(\cdot) \right)_{S^{n+1}}$ restricted to the support set $S^{n+1}$.

\newpage

\beginsupplement
\begin{figure}[h]
\begin{minipage}[b]{1.0\textwidth}
\centering
\includegraphics[width=0.85\linewidth]{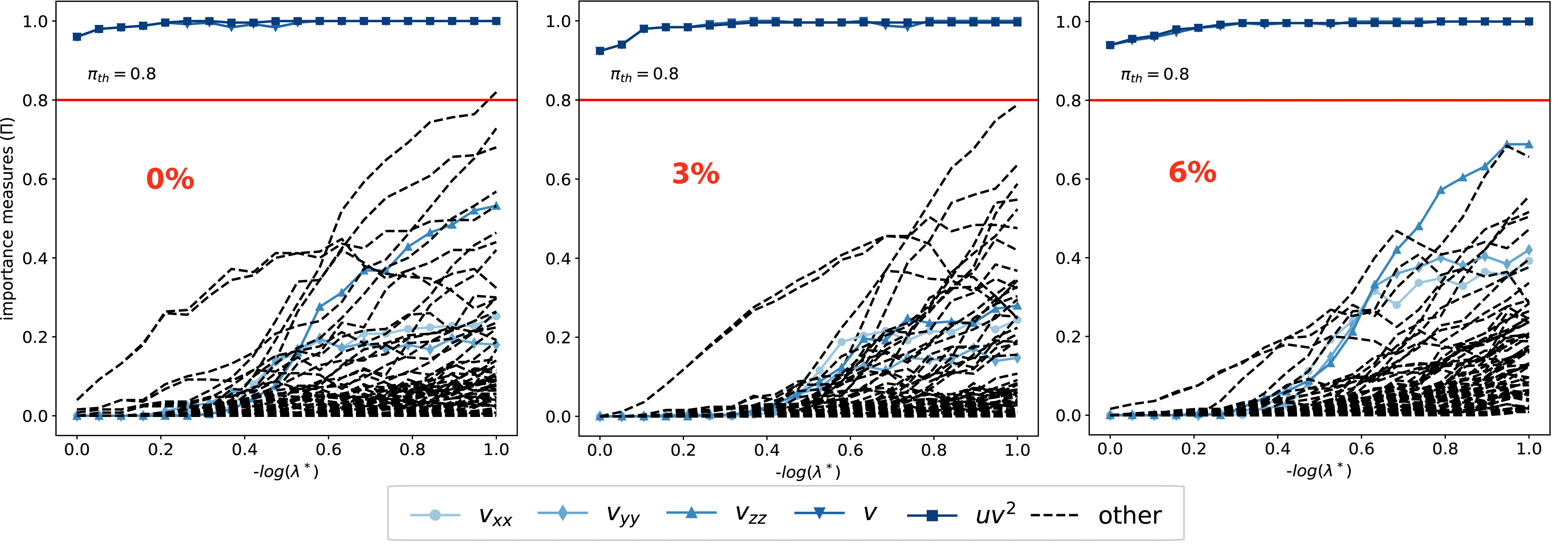}
\end{minipage}
\captionsetup{justification=raggedright,margin=0.5cm}
\caption{ \textbf{Model selection with PDE-STRIDE+STRidge for 3D Gray-scott $v-$component equation recovery :} The stability plots for the design $N=400, p=69$ show the separation of the true PDE components (in solid color) from the noisy components. The inference power of the PDE-STRIDE method is tested for additive Gaussian noise-levels $\sigma$ up-to $6\%$ (not shown). The PDE-STRIDE fails to identify the true support of the $v$-component PDE given the small diffusion coefficients $(1\textrm{E}^{-5})$ associated with the unidentified diffusive components. However, the plots show consistency for the reaction terms upto $6 \%$ additive Gaussian noise-levels. The inset at the bottom shows the colors correspondence with the PDE components.}
\label{Fig_gray_scott_v}
\end{figure}

\begin{table}[h]
\centering
\scalebox{0.9}{
\begin{tabular}{ |c|c|c|c| } 
 \hline
  & 1 & $ v (-0.067) $ & $ uv^2 (1.0) $ \\
 \hline
 clean & 0 & -0.0669  & 0.9999 \\ 
 \hline
 $2 \% $ & 0 & -0.0666 & 0.9910 \\ 
 \hline
 $4 \% $ & 0.0001  & -0.0667 & 0.9840\\
 \hline
  $6 \% $ & 0.0001  &  -0.0657 & 0.9677\\ 
 \hline
\end{tabular}}
\vspace{0.5em}
\captionsetup{justification=raggedright,margin=1.0cm}
\caption{ Coefficients of the recovered $v$-component Gray-Scott reaction diffusion equation for different noise levels. The stable components of the PDE from Figure \ref{Fig_gray_scott_v} are $\hat{S}_{stable} = \{ v ,  uv^2 \}$ }.
\label{coefficient_values_vGS}
\end{table}

\newpage

\subsection{Denoising technique}
For all the data-sets considered both from experiments and simulation, we use Singular Value decomposition (SVD) for de-noising. De-noising is achieved by identifying the \enquote{elbows} in singular values plots and applying a hardthreshold at the elbow to filter the noise \cite{donoho2013optimal, hansen1990truncated}.

\begin{figure}[h]
\centering
\includegraphics[width=0.6\linewidth]{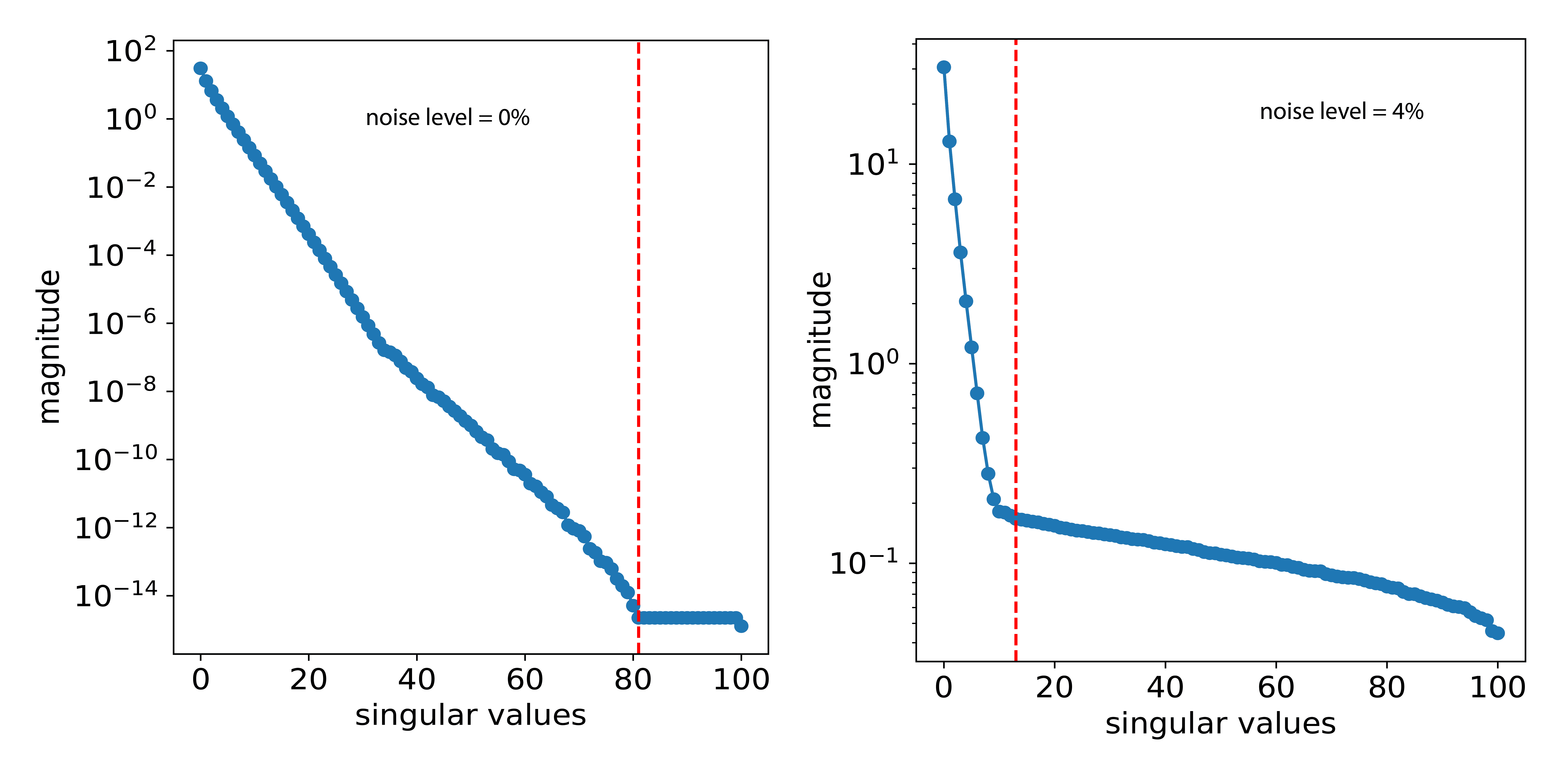}
\caption{ \textbf{ Truncated singular value decomposition (SVD) for denoising for 1D Burgers data-set }}
\label{Burgers_SVD}
\end{figure}

\begin{figure}[h]
\centering
\includegraphics[width=0.6\linewidth]{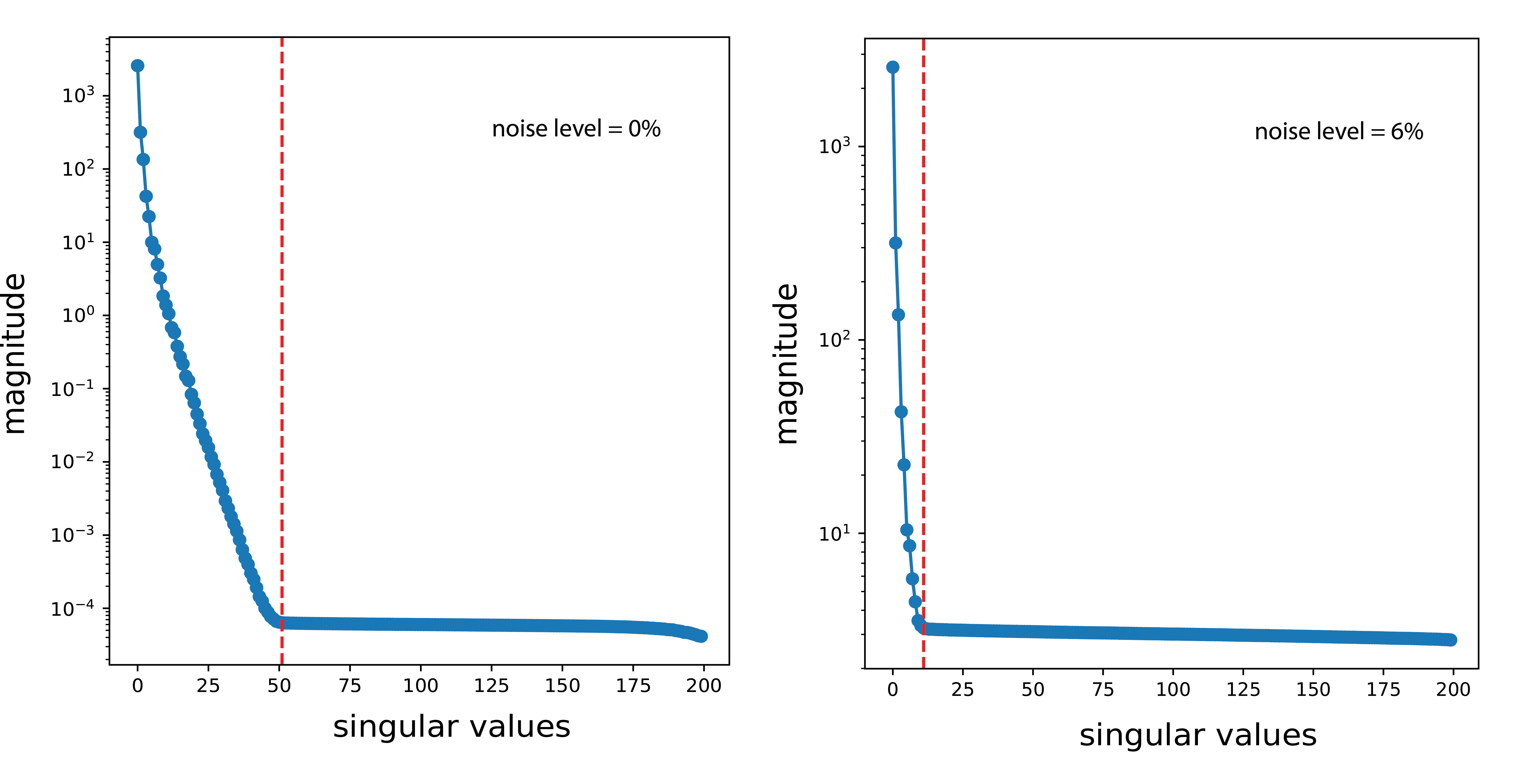}
\caption{ \textbf{ Truncated singular value decomposition (SVD) for denoising for 1D Burgers data-set }}
\label{GSRD_SVD}
\end{figure}

\newpage

\begin{figure}[h]
\begin{minipage}[b]{1.0\textwidth}
\centering
\includegraphics[width=0.85\linewidth]{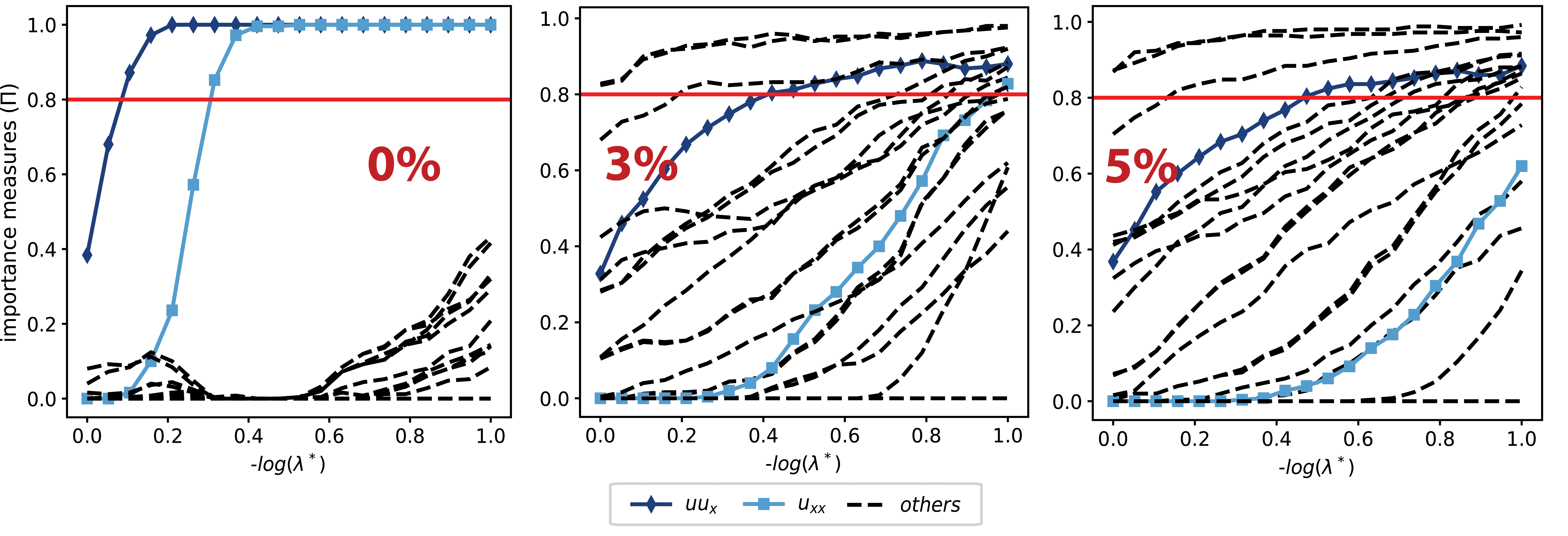}
\end{minipage}
\caption{ \textbf{Model selection with PDE-STRIDE+STRidge for 1D Burgers equation inference :}  \footnotesize The plots show the identification of the true PDE components (in solid color) from the noisy components. The Stability selection parameters are  $N=250, p=19$ with $B = 250$ repetitions. The statistical power of the algorithm is tested for additive Gaussian noise-levels upto $4\%$. All in the cases, perfect recovery was possible with the fixed threshold of $\pi^{th} = 0.8$ shown by the red solid line. The inset at the bottom shows the colors correspondence with the respective PDE components.}
\label{Burgers_results_STR_p20}
\end{figure}

\begin{figure}[h]
\begin{minipage}[b]{1.0\textwidth}
\centering
\includegraphics[width=0.85\linewidth]{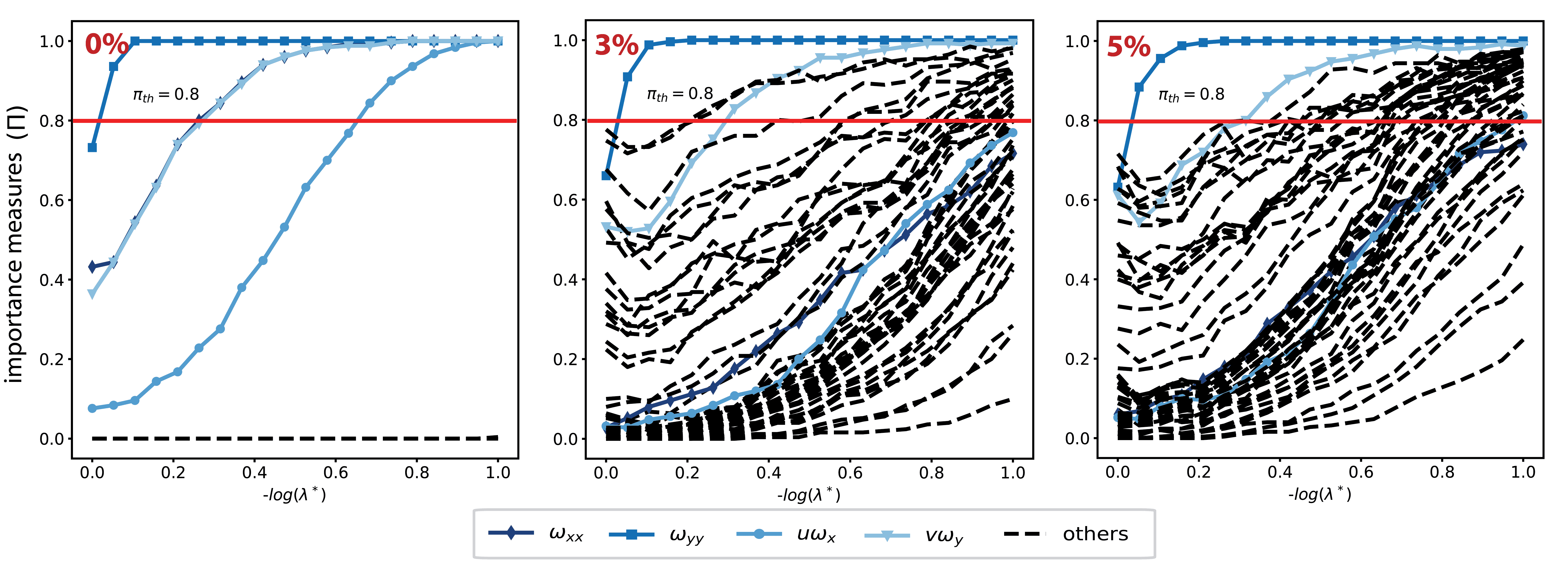}
\end{minipage}
\caption{ \textbf{Model selection with PDE-STRIDE+STRidge for 2D Vorticity transport equation inference :}  \footnotesize The plots show the identification of the true PDE components (in solid color) from the noisy components. The Stability selection parameters are  $N=500, p=48$ with $B = 250$ repetitions. The statistical power of the algorithm is tested for additive Gaussian noise-levels upto $5\%$. All in the cases, perfect recovery was possible with the fixed threshold of $\pi^{th} = 0.8$ shown by the dark blue solid line. The inset at the bottom shows the colors correspondence with the respective PDE components.}
\label{NS_STR_results_p48}
\end{figure}

\newpage

\begin{figure}[h]
\begin{minipage}[b]{1.0\textwidth}
\centering
\includegraphics[width=0.85\linewidth]{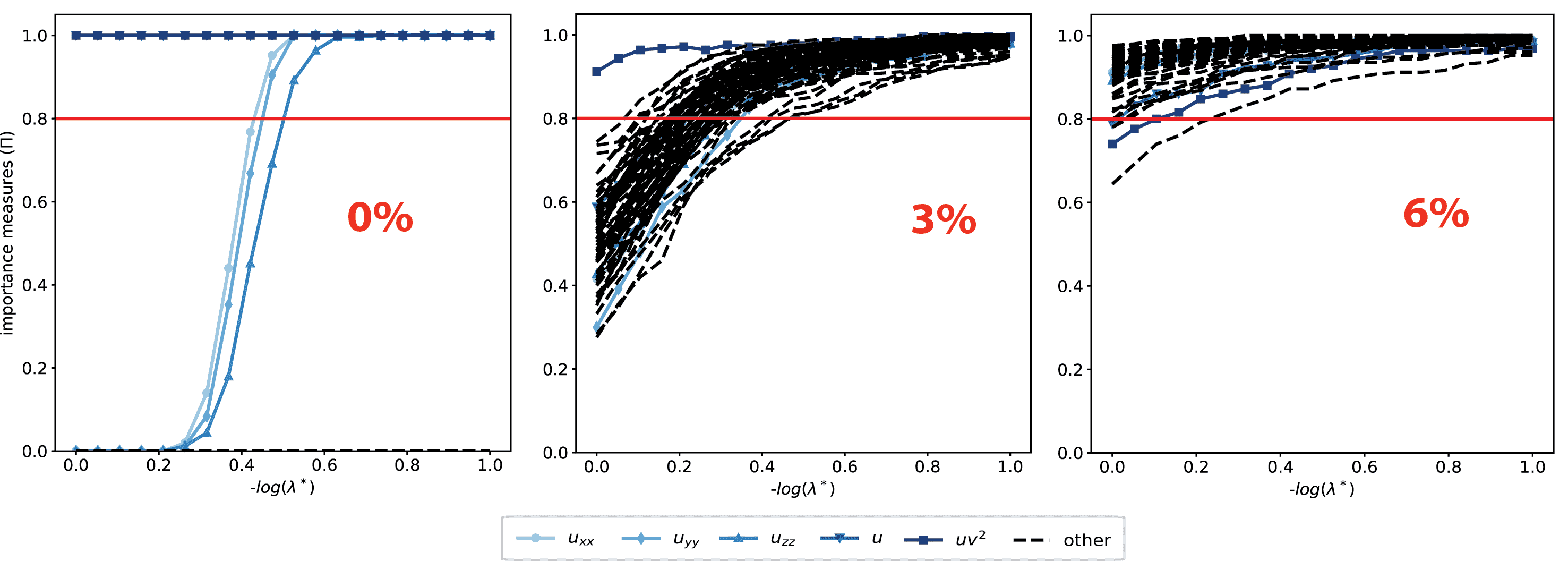}
\end{minipage}
\caption{ \textbf{Model selection with PDE-STRIDE+STRidge for inference of $u$-component of the Gray-Scott reaction diffusion equation :}  \footnotesize The plots show the identification of the true PDE components (in solid color) from the noisy components. The Stability selection parameters are  $N=400, p=69$ with $B = 250$ repetitions. The statistical power of the algorithm is tested for additive Gaussian noise-levels upto $6\%$. All in the cases, perfect recovery was possible with the fixed threshold of $\pi^{th} = 0.8$ shown by the dark blue solid line. The inset at the bottom shows the colors correspondence with the respective PDE components.}
\label{uGS_STR_results_p48}
\end{figure}

\begin{figure}[h]
\begin{minipage}[b]{1.0\textwidth}
\centering
\includegraphics[width=0.85\linewidth]{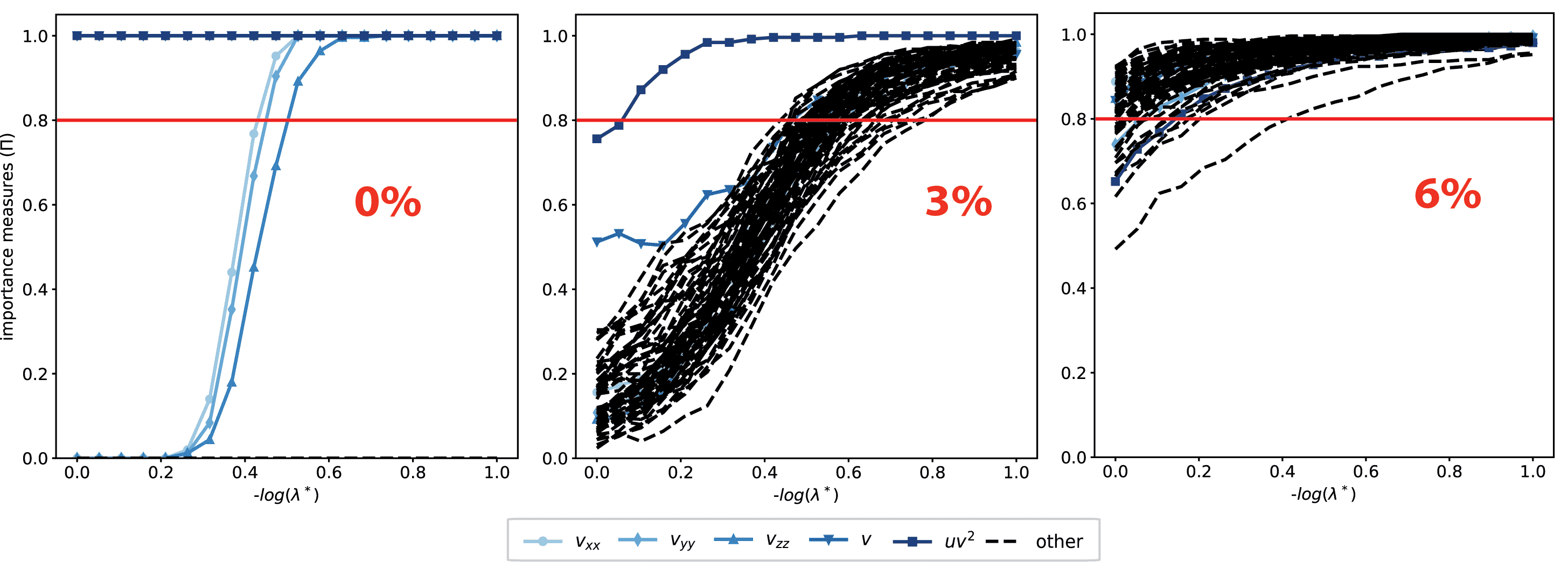}
\end{minipage}
\caption{ \textbf{Model selection with PDE-STRIDE+STRidge for inference of $v$-component of the Gray-Scott reaction diffusion equation :}  \footnotesize The plots show the identification of the true PDE components (in solid color) from the noisy components. The Stability selection parameters are  $N=400, p=69$ with $B = 250$ repetitions. The statistical power of the algorithm is tested for additive Gaussian noise-levels upto $6\%$. All in the cases, perfect recovery was possible with the fixed threshold of $\pi^{th} = 0.8$ shown by the dark blue solid line. The inset at the bottom shows the colors correspondence with the respective PDE components.}
\label{vGS_STR_results_p48}
\end{figure}

\newpage

\begin{figure}[h]
\begin{minipage}[b]{1.0\textwidth}
\centering
\includegraphics[width=0.85\linewidth]{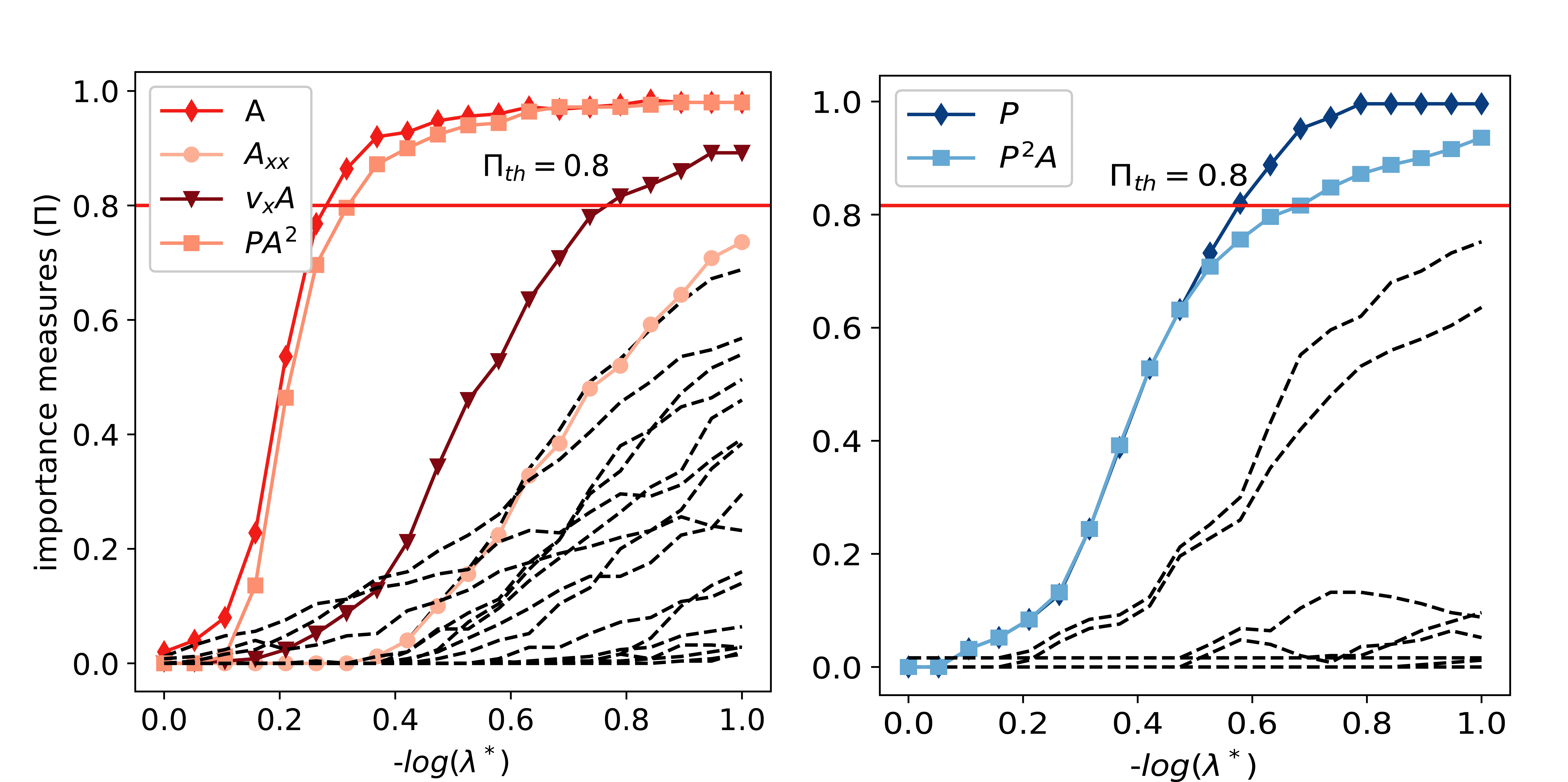}
\end{minipage}
\vspace{0.5em}
\caption{ \textbf{Model selection with PDE-STRIDE+IHT-d for PAR model inference :}  \footnotesize The stability results for the design $N=500,p=20$ for both aPAR (left) and pPAR (right) are shown. In the stability set for aPAR $\hat{S}_{stable} = \{ A, A^2P, \textrm{v}_{x} A \}$, we see advection dominant term $\textrm{v}_{x} A$ also appearing. The diffusion term $A_{xx}$ is also seen very close to the threshold $\pi_{th}=0.8$.}
\label{extra_PAR}
\end{figure}
\end{document}